\documentclass{article} 

\PassOptionsToPackage{round}{natbib}

\usepackage{authblk,breakcites}
\usepackage[round]{natbib}
\setcitestyle{authoryear}
\usepackage{microtype}
\usepackage{graphicx}
\usepackage{sidecap}
\usepackage{caption}
\usepackage{subfigure}
\usepackage{subcaption}
\usepackage{booktabs} 
\usepackage{multirow}
\usepackage[svgnames,dvipsnames]{xcolor}
\usepackage{float} 
\usepackage{wrapfig}
\usepackage{geometry}

\usepackage{hyperref}
\hypersetup{
    colorlinks=true,
    linkcolor=blue,
    filecolor=magenta,      
    urlcolor=cyan,
    citecolor=Emerald
    }

\usepackage{amsmath}
\usepackage{amssymb}
\usepackage{amsfonts}
\usepackage{mathtools}
\usepackage{amscd}
\usepackage{mathrsfs}
\usepackage{amsthm}
\usepackage[capitalize,noabbrev]{cleveref}
\usepackage{enumerate}

\usepackage{dsfont}

\usepackage{ulem}

\def \N{\mathbb{N}}
\def \R{\mathbb{R}}

\def \E{\mathbb{E}}
\def \F{\mathbb{F}}
\def \G{\mathbb{G}}

\def \a{\mathbf{a}}
\def \d{\mathbf{d}}

\def \x{\mathbf{x}}
\def \y{\mathbf{y}}

\def \P{\mathbb{P}}

\def \D{\mathbb{D}}

\def \1{\mathds{1}}
\def \0{\mathds{O}}

\def \Ac{\mathcal{A}}
\def \Bc{\mathcal{B}}
\def \Cc{\mathcal{C}}

\def \Ec{\mathcal{E}}
\def \Fc{\mathcal{F}}
\def \Gc{\mathcal{G}}

\def \Ic{\mathcal{I}}

\def \Lc{\mathcal{L}}
\def \Mc{\mathcal{M}}
\def \Nc{\mathcal{N}}

\def \Pc{\mathcal{P}}

\def \Rc{\mathcal{R}}
\def \Sc{\mathcal{S}}
\def \Tc{\mathcal{T}}

\def \Vc{\mathcal{V}}

\def \Pb{\textbf{P}}

\def \Db{\textbf{D}}

\def \Igen{{\Ic^{\rm gen}}}

\newcommand{\restr}[2]{#1_{\mkern 1mu \vrule height 2ex\mkern2mu #2}}

\def \eps{\varepsilon}

\def \eqsp{\;}

\def \l{\left}
\def \r{\right}

\providecommand{\keywords}[1]
{
  \small	
  \textbf{Keywords---} #1.
}
\providecommand{\MSC}[1]
{
  \small	
  \textbf{MSC Classification---} #1.
}

\theoremstyle{plain}
\newtheorem{theorem}{Theorem}[section]
\newtheorem{proposition}[theorem]{Proposition}
\newtheorem{lemma}[theorem]{Lemma}
\newtheorem{corollary}[theorem]{Corollary}
\theoremstyle{definition}
\newtheorem{definition}[theorem]{Definition}

\theoremstyle{remark}
\newtheorem{remark}[theorem]{Remark}

\newcounter{hypH}
\newenvironment{hypH}{
    \refstepcounter{hypH}
    \begin{itemize}
    \item[{\bf H\arabic{hypH}}]
    }
{\end{itemize}}

\newcommand{\eg}{\textit{e.g.}}
\newcommand{\ie}{\textit{i.e.}}

\def \CtrlStandard{\Sc}
\def \NatCtrl{\Rc^{\mathfrak{n}}}
\def \RelCtrl{\Rc^{\mathfrak{r}}}
\def \WeakCtrl{\Rc^{\mathfrak{w}}}
\def \CtrlRule{\Rc}
\def \CtrlRuleEps{\CtrlRule^\eps}
\def \SymmCtrl{\CtrlStandard^\mathfrak{s}}
\def \valueStrong{v^\mathfrak{s}}
\def \valueRelaxed{v^\mathfrak{r}}

\def \Leb{\Lambda}

\title{Controlled Interacting Branching Diffusion Processes: Relaxed Formulation in the Mean-Field Regime}

\author{Antonio Ocello}

\affil{{\small Centre de Recherche en Économie et de Statistiques (CREST), Groupe ENSAE-ENSAI, ENSAE Paris, Institut Polytechnique de Paris, 91120 Palaiseau, France}}

             \date{\today}

\usepackage[colorinlistoftodos]{todonotes} 

\begin{document}

\maketitle

\begin{abstract}
    The focus of this article is studying an optimal control problem for branching diffusion processes. Initially, we introduce the problem in its strong formulation and expand it to include linearly growing drifts. Then, we present a relaxed formulation that provides a suitable characterization based on martingale measures. Considering weak controls, we prove they are equivalent to strong controls in the relaxed setting, and establish the equivalence between the strong and relaxed problem, under a Filippov--type convexity condition. Furthermore, by defining control rules, we can restate the problem as the minimization of a lower semi-continuous function over a compact set, leading to the existence of optimal controls both for the relaxed problem and the strong one. Finally, with a useful embedding technique, we show that the value function solves a system of HJB equations, establishing a verification theorem. We then apply it to a linear-quadratic example and a kinetic one.
\end{abstract}

\noindent \MSC{93E20, 60J60, 60J80, 35K10, 60J70, 60J85}

\noindent \keywords{Stochastic control, relaxed control, branching diffusion processes, martingale representation}


\section{Introduction}

    The focus of this paper is on populations that are optimally controlled under a mean-field constraint. Specifically, we aim to show the presence of a strong control for controlled branching diffusions and to describe the optimal dynamics in this regime. Branching diffusion processes, originally introduced in \citet{Sk64, INW69, INW681, INW682}, describe the evolution of particles whose spatial motion is governed by stochastic differential equations. They have since been extensively studied, notably for their role in the probabilistic representation of semilinear PDEs \citep[see, \eg,][]{HLTT} and in applications such as regularized unbalanced optimal transport \citep{baradat2021regularized}. More broadly, these processes belong to the class of measure-valued processes, which have been a central focus of research since the late 20th century. Seminal works such as \citet{Roelly-Meleard:discont_measure_val_branching, Roelly:Rouault:Branchements_Spatiaux, Roelly:convergence_measure_val_processes} introduced these processes as solutions to martingale problems, a perspective formalized in the diffusion case by \citet{ethier2009markov}. This martingale approach, further developed in \citet{Dawson} to encompass Fleming–Viot processes and superprocesses, provides both an abstract and powerful framework together with convergence criteria that will play a key role in this paper.

Several examples of \emph{optimal control for branching processes} are discussed in the literature \citep{Ustunel,Nisio,claisse18,kharroubi2024stochastic,kharroubi2024optimal, claisse2023mean,bethencourt2025brownian}. 
The initial reference to these objects can be found in \citet{Ustunel}, wherein their modelling employs a disjoint topological sum of Euclidean space. The control, living within a compact space, solely affects the drift of spatial movement. The author permits each particle to potentially be influenced by any other living particle, without imposing any additional assumptions on the structure of these interactions. Moreover, the running cost yields a high degree of generality as well, leading to a correspondingly complex differential characterization.
By selecting the cost function as the product of functions associated with the living particles at the terminal time, \citet{Nisio} employs controlled branching processes as a probabilistic tool to examine a specific group of parabolic Bellman equations. In this study, the control, still confined to a compact set, influences both drift and volatility. A Hamilton--Jacobi--Bellmann (HJB) equation is identified, establishing that the value function represents its unique (viscosity) solution.

In \citet{claisse18}, the author goes further in the analysis of this setting. Initially, the controlled processes are described as measure-valued processes. Using Ulam--Harris--Neveu labelling \citep[see, \eg,][]{Ulam-Harris_example} to describe the genealogy of the particles, the author introduces a label set that assists in defining the branching events. A set of Brownian motions and Poisson random measures, indexed by these labels, are utilized to provide a strong formulation for the controlled branching processes. This facilitates proving the well-posedness for dynamics where drift, volatility, branching rate, and branching mechanisms are not only controlled but also dependent on the position of each particle.
While these coefficients are still assumed to be bounded, the control space is no longer necessarily compact. Since the dynamics are coupled only through the control, the product structure of the cost yields a branching property that converts the problem into a finite-dimensional one. A PDE characterization of the value function is then obtained, leveraging the differential properties of the Euclidean space where each single particle is defined. In \citet{kharroubi2024stochastic}, a similar approach is also employed. Here, the symmetry of the reward function is again used to establish a different branching property that allows for finite-dimensional rewriting.

In a companion paper \citep{ocello2025controlled1}, we addressed the problem of optimal control for branching diffusion processes under general interaction assumptions. This was carried out in full generality using HJB techniques, which characterize the value function through the framework of viscosity solutions. Moreover, Proposition 4.2 in that work demonstrates that, under a \emph{mean-field interaction} assumption, it is optimal to restrict attention to symmetric controls.
The mean-field setting is built on two fundamental assumptions. The first is \emph{anonymity}, which requires that agents react only to the empirical distribution of the population’s positions, without distinguishing between individual identities. The second is \emph{homogeneity}, which stipulates that an agent’s behavior is independent of its specific label. These assumptions naturally restrict the generality of the model parameters: anonymity is reflected through dependence on the empirical distribution $\pi(\lambda)$ instead of the full configuration $\lambda$, while homogeneity is enforced by requiring invariance with respect to the particle index. 

In the present work, we aim to complete the analysis of the mean-field regime by addressing the existence of optimal controls through the lens of a \emph{relaxed formulation}. 
This approach enables us to overcome compactness issues inherent in the space of strict controls and to establish the existence of minimizers.
Similar methodology has been used in mean-field control theory \citep{lacker2017limit,Bahlali-et-al:Existence_opt_MFeq} or branching populations dynamics \citep{claisse2023mean}. Our approach differs from \citet{claisse2023mean} as they make large use of the indexation with respect to the label set.
The relaxed control formulation shows key links with the martingale problem viewpoint of the control problem and exhibit strong theoretical properties, as emphasized in \citet{ElK:Tan:Capacities2}.
We follow the approach of \citet{elkaroui1987compactification} and \citet{haussmann1990existence}, which involves a relaxed formulation of the problem. In this setting, the control process is represented as a probability measure on the action space, depending on both time and space, and satisfying specific well-posedness constraints. This formalism introduces different descriptions of the control problem, namely control rules and natural controls, allowing for greater flexibility and easier manipulation of the controlled dynamics. Proving that a control rule (resp. natural control) with a lower cost can be constructed from any relaxed control (resp. control rule), we establish the equivalence between strong and relaxed problems. This is achieved through a Filippov-type convexity condition \citep{Filippov}, showing that any relaxed control can be associated with a weak control of lower cost, allowing us to restrict the search for minimizers to a compact set under a suitable topology.
Therefore, we prove that the cost functional is lower semicontinuous in the control rule framework. 
This rewrites the original optimization as the minimization of a lower semicontinuous function over a compact space, establishing the existence of optimal values and controls.


This article lays the groundwork for applying a broad class of numerical methods, as the relaxed formulation provides the rigorous foundation many of them rely on. In particular, in Markov chain approximation methods developed by Kushner and Dupuis \citep{kushner1990numerical}, convergence to the true control problem relies on the compactness of the relaxed control space. Similarly, in the context of numerical resolution of HJB equations via finite difference or finite element schemes, the compactness and convexity of the relaxed formulation guarantee the existence of minimizers of the Hamiltonian and ensure the stability and convergence of the approximation \citep{barles1991convergence, krylov87}. Moreover, in stochastic approximation and reinforcement learning, relaxed controls naturally emerge as the limiting objects of randomized or exploratory strategies, furnishing the right topological framework to establish convergence results \citep[see, e.g.,][]{kushner2003stochastic, borkar2008stochastic}. These developments, now accessible in the context of controlled branching processes, offer a promising path forward for both theoretical analysis and computational practice.
Beyond numerical aspects, the martingale problem point of view also lays the foundation for studying scaling limits of controlled branching processes, with one such limit rigorously developed in \citet{ocello2025controlled}. These developments, now accessible in the context of controlled branching processes, offer a promising path forward for both theoretical analysis and computational practice.

The remainder of the paper is organized as follows. In \Cref{section:setting}, we introduce the model and assumptions, together with the strong formulation of controlled branching processes and the associated control problem in the mean-field regime. \Cref{section:relax_form} develops the relaxed formulation of this problem, through the martingale problem perspective, and introduces the natural and weak control classes. Under a Filippov-type convexity condition, we show how any natural control can be associated with a weak one of smaller cost. Finally, in \Cref{section:Existence_optimum}, we establish the equivalence between relaxed and strong formulations, introduce the notion of control rules, and prove that the original control problem reduces to optimizing a lower semicontinuous functional over a compact set—thus ensuring the existence of a minimizer.

\section{Setting}
\label{section:setting}

\subsection{Notation}

\paragraph{Finite measures.}
For a Polish space $(\Ec,d)$ with $\Bc(\Ec)$ its Borelian $\sigma$-field, we write $C_b(\Ec)$ (resp. $C_0(\Ec)$) for the subset of the continuous functions that are bounded (resp. that vanish at infinity), and $M(\Ec)$ (resp. $\Pc(\Ec)$) for the set of Borel positive finite measures (resp. probability measures) on $\Ec$. We equip $M(\Ec)$ with weak* topology, \ie, the weakest topology that makes continuous the maps $M(\Ec)\ni\lambda\mapsto\int_\Ec \varphi(x) \lambda(dx)$, for $\varphi\in C_b(\Ec)$. We denote $\langle \varphi,\lambda\rangle = \int_\Ec \varphi(x) \lambda(dx)$, for $\lambda\in M(\Ec)$ and $\varphi\in C_b(\Ec)$.

Denote also by $M^1(\Ec)$ the subspace of measures with finite first order moment, \ie, the collection of all $\lambda\in M(\Ec)$ such that $\int_\Ec d(x,x_0)\lambda(dx)<\infty$, for some $x_0\in \Ec$. The weak* topology can be metrized in $M^1(\Ec)$ by the Wasserstein type metric $\d_{1,\Ec}$, as introduced in Appendix B of \citet{claisse2023mean}. This means that, if $\partial$ is a cemetery point, we consider first $\bar \Ec$ the enlarged space $\bar \Ec := \Ec\cup \{\partial\}$. Defining $d(x,\partial):=d(x,x_0)+1$, we have that $(\bar \Ec,d)$ is Polish. For $m\in\R_+$, we consider the Wasserstein distance $\d_{1,\Ec,m}$, on the space $M^1_m(\bar \Ec)$ defined as
\begin{align*}
    M^1_m(\bar \Ec) := \{\lambda \in M^1(\bar \Ec) :  \lambda(\bar \Ec) = m\}\eqsp,
\end{align*}
as follows
\begin{align*}
    \d_{1,\Ec,m}(\lambda,\lambda^\prime) = \inf_{\pi\in\Pi(\lambda,\lambda^\prime)}\int_{\bar \Ec \times \bar \Ec}d(x,y)\pi(dx,dy)
    \eqsp, \quad \text{ for }\lambda,\lambda^\prime\in M^1_m(\bar \Ec)\eqsp,
\end{align*}
with $\Pi(\lambda,\lambda^\prime)$ the collection of all non-negative measures on $\bar \Ec \times \bar \Ec$ with marginals $\lambda$ and $\lambda^\prime$.
The distance $\d_{1,\Ec}$ on $M^1(\Ec)$ is now defined as
\begin{align*}
    \d_{1,\Ec}(\lambda,\lambda^\prime) = \d_{1,\Ec,m}\left(\bar\lambda_m,\bar\lambda^\prime_m\right)\eqsp, \quad \text{ for }\lambda,\lambda^\prime\in M^1_m(\Ec)\eqsp,
\end{align*}
with $m \geq \lambda(\Ec)\vee \lambda^\prime(\Ec)$, $\bar\lambda_m(\cdot):=\lambda(\cdot\cap \Ec) + (m-\lambda(\Ec))\delta_\partial(\cdot)$,and $\bar\lambda^\prime_m(\cdot):=\lambda^\prime(\cdot\cap \Ec) + (m-\lambda^\prime(\Ec))\delta_\partial(\cdot)$.
As proven in Lemma B.1 of \citet{claisse2023mean}, this definition does not depend on the choice of $m$. Moreover, for some $x_0\in \Ec$, we have the natural bound
\begin{align}\label{eq:bound_d_p_E}
    \d_{1,\Ec}(\lambda,\delta_{x_0})\leq \int_\Ec d(x,x_0)\lambda(dx) + \langle1,\lambda\rangle
    \eqsp,
    \quad \text{ for }\lambda\in M^1(\Ec)\eqsp.
\end{align}

Finally, we write $\Nc(\Ec)$, for the space of weak measures in $\Ec$, \ie,
\begin{align*}
    \Nc(\Ec) := \left\{\sum_{i=1}^m\delta_{x_i} ~:~m\in\N, x_i\in \Ec \text{ for }i\leq m\right\}\eqsp,
\end{align*}
a weakly* closed subset of $M(\Ec)$.

\paragraph{Label set.} We use Ulam--Harris--Neveu labelling to consider the genealogy of the particles.
Consider the set of labels 
\begin{align*}
    \Ic := \{\varnothing\}\cup\bigcup_{n=1}^{+\infty}\N^n\eqsp.
\end{align*}
Denote by $\varnothing$ the \textit{mother particle}, and $i=i_1\cdots i_n$ the multi-integer $i=(i_1,\ldots,i_n)\in\N^n$, $n\geq 1$.  For $i=i_1\cdots i_n\in \N^n$ and $j=j_1\cdots j_m\in \N^m$, we define their concatenation is $i j\in\N^{n+m}$ by $i j =i_1\cdots i_n j_1\cdots j_m$, and extend it to the entire $\Ic$ by $\varnothing i = i \varnothing=i$, for all $i\in\Ic$. When a particle $i=i_1\cdots i_n\in \N^n$ gives birth to $k$ particles, the off-springs are labelled $i0,\ldots,i(k-1)$. Moreover, if $\Vc\subset\Ic$ was the set of alive particles, after the branching event on the branch $i\in\Vc$, we have that the new set of alive particles become $\Vc^i_k$, with
\begin{align}
\label{eq:def:V_i_k}
    \Vc^i_k := \Vc\setminus\{i\}\cup \l\{i0,\dots,i(k-1)\r\}\eqsp.
\end{align}

Consider the partial ordering $\preceq$ (resp. $\prec$) by
\begin{align*}
    i\preceq j ~ \Leftrightarrow ~ \exists \ell\in\Ic~:~j=i\ell
    \qquad
    \left(\textrm{resp.}~i\prec j ~ \Leftrightarrow ~\exists \ell\in\Ic\setminus \{\varnothing\}~:~j=i\ell\right)\eqsp,
\end{align*}
for $i,j\in\Ic$.
We endow $\Ic$ with the discrete topology, generated by the distance
\begin{align*}
d^\Ic(i,j) := \sum_{\ell = p+1}^n (i_\ell +1) + \sum_{\ell' = p+1}^m (j_{\ell'} +1)\;,
\qquad\text{ for }
i=i_1\cdots i_{n}\in\N^n, \;j=j_1\cdots j_{m} \in \N^m\eqsp,
\end{align*}
with $p = \max\{\ell\geq1:i_\ell=j_\ell\}$ the generation of the greatest common ancestor. Denote $i\wedge j= i_0\cdots i_p$ and write $|i|:= d^\Ic(i,\varnothing)$, for $i\in\Ic$. Moreover, define the total ordering $\leq$ on $\Ic$ as $i\leq j$ if $i\preceq j$ or $i_{p+1}<j_{p+1}$.


From the definition of $\Vc^i_k$, note that not all possible combinations of indeces are considered when describing a population. Let $\Ic^{gen}$ be the space of admissible configurations of indeces, defined as
\begin{align*}
    \Ic^{gen}
    :=&\Big\{
        \Vc\eqsp:\eqsp \Vc\subseteq \Ic \text{ finite},\eqsp i\nprec  j,\text{ for }i,j\in\Vc
    \Big\}\eqsp.
\end{align*}
As $\Igen$ is a subset of $\mathscr{P}_{\mathrm{fin}}(\Ic)$ the set of all finite subsets of $\Ic$, it is a countable set.
For $\Vc\in\Igen$, denote $\mathfrak{S}_\Vc$, the set of permutations of $\Ic$ that send $\Vc$ to an admissible configuration in $\Igen$, defined as
\begin{align*}
    \mathfrak{S}_\Vc := 
    \bigl\{
        \mathfrak{s}\in \mathrm{Sym}(\Ic)\;:\;\mathfrak{s}\cdot \Vc \in \Igen
    \bigr\}\eqsp,
\end{align*}
where $\mathrm{Sym}(\Ic)$ is the permutation group of $\Ic$ and the action on subsets is $\mathfrak{s}\cdot \Vc \;:=\; \{\mathfrak{s}(i): i\in \Vc\}\subset \Ic$. Moreover, denote $\mathfrak{s}\cdot\check\lambda$ to be $\mathfrak{s}\cdot\check\lambda := \sum_{i\in\Vc}\delta_{(\mathfrak{s}(i),x_i)}$, for $\check\lambda = \sum_{i\in\Vc}\delta_{(i,x_i)}\in E$.

\paragraph{State and control space.} Take $E\subset \Nc(\Ic\times\R^d)$ as
\begin{align*}
    E:=&\l\{
        \sum_{i\in \Vc}\delta_{(i,x_i)}\eqsp : \eqsp \Vc\in \Ic^{gen},\; x_i\in\R^d
    \r\}\eqsp.
\end{align*}
Note that $\Nc(\R^d)$ is a closed set of $M^1(\R^d)$ with respect to the distance $\d_{1,\R^d}$. This is due to the fact that $\Nc(\R^d)$ is weakly*-closed and, from Lemma B.2 in \citet{claisse2023mean}, convergence in $M^1(\R^d)$ entails weak*-convergence to some $\lambda\in \Nc(\R^d)\subseteq M^1(\R^d)$. Therefore, combining this with the fact the $E$ is weakly*-closed \citep[see, $e.g$, Proposition A.7,][]{kharroubi2024stochastic} and $\Ic$ is equipped with discrete topology, we have that $E$ is also a closed set of $M^1(\Ic\times \R^d)$.

Define now this projection map $\pi:E\to \Nc(\R^d)$ 
as
\begin{align*}
    \pi:E\ni \sum_{i\in \Vc}\delta_{(i,x_i)}\mapsto \sum_{i\in \Vc}\delta_{x_i}
    \eqsp.
\end{align*}
Fix $\check\lambda=\sum_{i\in \Vc}\delta_{(i,x_i)}$, $\check\lambda^\prime=\sum_{i\in \Vc}\delta_{(i,y_i)}\in E$. Using the characterisation of the distance $\d_{1,\Ic\times \R^d}$ of Lemma B.1 in \citet{claisse2023mean}, we obtain
\begin{align}
\label{eq:bound_d_1_Rd-measures-to_R_d}
    \d_{1,\Ic\times \R^d}\left(\check\lambda,\check\lambda^\prime\right) = 
    \sup_{\varphi\in\text{Lip}^0_1(\Ic\times \R^d)}\sum_{i\in\Vc}\left|\varphi(i,x_i) - \varphi(i,y_i) \right| \leq \sum_{i\in\Vc}|x_i-y_i| = \|\Vec{x}_\Vc-\Vec{y}_\Vc\|_{1,d|\Vc|}
    \eqsp,
\end{align}
where $\Vec{x}_\Vc=(x_i)_{i\in \Vc}$ is the vector of $\R^{d|\Vc|}$ taken in the order induced by the total ordering $\leq$ on $\Ic$,
$\text{Lip}^0_1(\R^d)$ denote the collection of all functions $\varphi : \Ic\times \R^d\to \R$ with Lipschitz constant smaller or equal to $1$ and such that $\varphi(0) = 0$ and $\|\cdot\|_{1,n}$ denotes the $L^1$-distance in $\R^{n}$, for $n\in\N$.
Using Cauchy--Schwarz inequality, we can also bound the distance $\d_{1,\Ic\times \R^d}$ by
\begin{align}
\label{eq:bound_d_1_Rd-measures-to_R_d-2}
    \d_{1,\R^d}\left(\check\lambda,\check\lambda^\prime\right) \leq \sqrt{|\Vc|}\eqsp \|\Vec{x}_\Vc-\Vec{y}_\Vc\|_{2,d|\Vc|}
    \eqsp,
\end{align}
where $\|\cdot\|_{2,n}$ denotes the $L^2$-distance in $\R^{n}$, for $n\in\N$.

Let $\Tc_{t,s}$ denotes the
collection of all stopping times valued in $[t,s]$.
Fix $m\in\R^d$. Take a closed subset $A$ of $\R^m$ representing the set of actions.

\paragraph{Càdlàg paths.}
Denote by $\Db^d = \D([t, +\infty), M^1(\R^d))$ (resp. $\D([0,T];E)$) the space of càdlàg functions from $[t, +\infty)$ to $M^1(\R^d)$ (resp. $E$), equipped with the Skorokhod topology $d_{\Db^d}$  (resp. $d_{E}$) associated with the metric $\d_{1,\R^d}$ (resp. $\d_{1,\Ic\times \R^d}$), which makes it complete \citep[see, \eg,][]{billingsley2013convergence}.
    
\subsection{Branching diffusion processes in the Mean-Field regime}
\label{Section:strong_form}
Fix a finite time horizon $T > 0$. 
Let $(\Omega,\Fc,\P)$ be a probability space supporting two independent families $\{W^i\}_{i\in \Ic}$ and $\{Q^i\}_{i\in \Ic}$ of mutually independent processes. Let $W^i$ be a $d^\prime$-dimensional Wiener processes, and $Q^i(ds dz)$ a Poisson random measure on $[0,T]\times\R_+$ with intensity measure $dsdz$.
Let $\F=\{\Fc_t\}_{t\geq 0}$ be the filtration generated by these processes, \ie, the (right-continuous) completion of the $\sigma$-algebra $\G=\{\Gc_t\}_{t\geq 0}$ with
\begin{align*}
    \Gc_t := \sigma\l(W^ i_s, Q^i([0,s]\times C)~:~ s \leq t,~ i\in \Ic,~ C\in \Bc(\R_+)\r)
    \eqsp.
\end{align*}
Moreover, let $\Fc_\infty$ (resp. $\Gc_\infty$) be the $\sigma$-algebra generated by $\bigcup_{t\geq 0}\Fc_t$ (resp. $\bigcup_{t\geq 0}\Gc_t$).

Consider the following parameter of models
\begin{align*}
(b,\sigma, \gamma, p_k
)\eqsp :\eqsp  \R^d\times \Nc(\R^d) \times A \to \R^d \times \R^{d\times d^\prime}\times \R_+ \times [0,1]\eqsp,
\end{align*}
for $k\geq0$, such that $\sum_{k \geq 0}p_k(x,\lambda,a)=1$, for $(x,\lambda,a)\in \R^d\times\Nc(\R^d)\times A$.
Let $\Phi$ be the generating function of $(p_k)_k$, \ie,
\begin{align*}
    \Phi(s,x,\lambda,a) = \sum_{k=0}^\infty p_k(x,\lambda,a) s^k\eqsp, \quad \text{ for }(s,x,\lambda,a)\in [0,1]\times \R^d\times\Nc(\R^d)\times A\eqsp.
\end{align*}
We now introduce the following assumptions on these parameters.

\begin{hypH}
\label{hypH:model_parameters}
    \begin{enumerate}[(i)]
        \item Suppose that $b$ and $\sigma$ are Lipschitz continuous in $(x,\lambda)$ uniformly in $a$, \ie, there exists $L>0$ such that
        \begin{align}
        \label{eq:bound_b_sigma_Lipschitz}
            \left|b(x,\lambda,a)-b(x',\lambda^\prime,a)\right|+\left|\sigma(x,\lambda,a)-\sigma(x',\lambda^\prime,a)\right|\leq
        L (\|x-x'\|_{2,d}+\d_{1,\R^d}(\lambda,\lambda^\prime))\eqsp,
        \end{align}
        for $x,x'\in\R^d$, $\lambda,\lambda^\prime\in \Nc(\R^d)$, and $a\in A,i\in\Ic$.
        \item Suppose that $\sigma$ and $\gamma$ are uniformly bounded, and $b$ has linear growth in $(x,a)$ while bounded in $(\lambda,i)$, \ie, there exists $C_\sigma, C_\gamma, C_b >0$ such that
        \begin{align}\label{eq:bound_b_sigma_gamma}
            \left|b(x,\lambda,a)\right|\leq C_b (1+ |x| + |a|)\eqsp,\qquad
            \left|\sigma(x,\lambda,a)\right|\leq C_\sigma\eqsp,\qquad
            \gamma(x,\lambda,a)\leq C_\gamma\eqsp,
        \end{align}
        for $(x,\lambda,a)\in \R^d\times\Nc(\R^d) \times A$.
        \item Suppose that the first and second order moments related to $(p_k)_k$ are uniformly bounded, \ie, there exist two constants $C^1_\Phi, C^2_\Phi>0$ such that
        \begin{align}
        \label{eq:bound:order1_2_Phi}
            \begin{split}
                \partial_s \Phi(1,x,\lambda,a)&=\sum_{k\geq1}k p_k(x,\lambda,a)\leq C^1_\Phi\eqsp,
                \qquad
                \partial^2_{ss} \Phi(1,x,\lambda,a)=\sum_{k\geq1}k(k-1) p_k(x,\lambda,a) \leq C^2_\Phi\eqsp,
            \end{split}
        \end{align}
        for $(x,\lambda,a)\in \R^d\times \Nc(\R^d)\times A$.
    \end{enumerate}
\end{hypH}

Consider the following definition of \emph{standard strong control}, which coincides with the notion of admissible control introduced in \citet{ocello2025controlled1}.

\begin{definition}[Standard strong control]
    We say that $\beta = (\beta^i)_{i\in \Ic}$ is a \emph{standard strong control}, and we denote $\beta\in\CtrlStandard$, if $\beta$ is an $\G$-predictable process valued in $A^\Ic$, such that
    \begin{align}\label{eq:bound_sup_beta_2}
        \E^{\P}\left[\int_t^T \sup_{i\in\Ic} |\beta_s^i|^2 ds\right]<\infty\eqsp.
    \end{align}
\end{definition}

Fix and initial condition $(t,\check\lambda)\in[0,T]\times E$ and  a standard control $\beta = (\beta^i)_{i\in\Ic}\in\CtrlStandard$. We describe the \textit{controlled branching diffusion} $\xi^{t,\check\lambda;\beta}$ as the measure-valued process
\begin{align*}
    \xi^{t,\check\lambda;\beta}_s = \sum_{i\in \Vc^{t,\check\lambda;\beta}_s}\delta_{(i,Y^{i,\beta}_s)}\eqsp,
\end{align*}
where $Y^{i,\beta}_s$ is the position of the member with label $i\in\Ic$, and $\Vc^{t,\check\lambda;\beta}_s$ the set of alive particles at time $s$.

Let $L$ be the generator (associated with the spatial motion of each particle) defined on $\varphi\in C^2_b(\R^d)$ as
\begin{align*}
    L \varphi (x,\lambda,a) = b(x,\lambda,a)^\top D \varphi(x) + \frac{1}{2}\text{Tr}\left(
    \sigma\sigma^\top(x,\lambda,a)D^2 \varphi(x)
    \right)\eqsp, \qquad\text{ for }(x,\lambda,a)\in \R^d\times \Nc(\R^d)\times A\eqsp,
\end{align*}
where $D\varphi$ and $D^2\varphi$ denote gradient and Hessian of the function $\varphi$.
A representation of branching process in the MF-regime is given by the following SDE 
\begin{align}
\label{SDE:strong}
    \begin{split}
        \langle\varphi,\xi^{t,\check\lambda;\beta}_s\rangle =&~ \langle\varphi,\lambda\rangle+
        \int_t^s\sum_{i\in \Vc^{t,\check\lambda;\beta}_u}D \varphi(i,Y^{i,\beta}_u)^\top\sigma\left(Y^{i,\beta}_u,\pi(\xi^{t,\check\lambda;\beta}_u),\beta^i_u\right)dB^i_u
        +
        \int_t^s\sum_{i\in \Vc^{t,\check\lambda;\beta}_u}L \varphi\left(Y^{i,\beta}_u,\pi(\xi^{t,\check\lambda;\beta}_u),\beta^i_u\right)du \\
        & + \int_{(t,s]\times\R_+}\sum_{i\in \Vc^{t,\check\lambda;\beta}_{u-}}\sum_{k\geq0}
        (k-1)\varphi(i,Y^{i,\beta}_{u-}) \1_{I_k\left(Y^{i,\beta}_{u-},\pi(\xi^{t,\check\lambda;\beta}_{u-}),\beta^i_{u}\right)}(z)
        Q^i(dudz)\eqsp,
    \end{split}
\end{align}
with
\begin{align*}
    I_k(x,\pi(\check\lambda),a) = \Bigg[\gamma(x,\pi(\check\lambda),a)\sum_{\ell=0}^{k-1}p_\ell(x,\pi(\check\lambda),a),\gamma(x,\pi(\check\lambda),a)\sum_{\ell=0}^{k}p_\ell(x,\pi(\check\lambda),a)\Bigg)\eqsp,
\end{align*}
for all $(x,\check\lambda,a)\in \R^d\times E \times A$, $k\in\N$, with the value of an empty sum being zero by convention. Notice that $(I_k(x,\pi(\check\lambda) ,a))_{k\in \N}$ forms a partition of the interval $[0, \gamma(x,\pi(\check\lambda),a))$.

Under Assumption H\ref{hypH:model_parameters}, we inherit the existence of controlled branching diffusions satisfying \eqref{SDE:strong} from Proposition 2.2 in \citet{ocello2025controlled1}, which extends the result of Proposition 2.1 in \citet{claisse18}. This ensures both the well-posedness of the dynamics under standard strong controls and the validity of the moment bounds required for the control problem.

\begin{proposition}
\label{prop:existence_strong_branching}
    Let $(t,\check\lambda) \in [0,T]\times E$ and $\beta\in\CtrlStandard$. Suppose Assumption H\ref{hypH:model_parameters} holds. Then, there exists a unique (up to indistinguishability) càdlàg and adapted process $(\xi^{t,\check\lambda;\beta}_s)_{s\geq t}$ satisfying~\eqref{SDE:strong} such that $\xi^{t,\check\lambda;\beta}_t=\check\lambda$. In addition, we have
    \begin{align}
        \label{eq:non-explosion-moment1_mass}
        \E^{\P}\left[\sup_{u\in[t,t+h]}|\Vc^{t,\check\lambda;\beta}_u|\right]\leq&~ \langle1,\check\lambda\rangle~ e^{C_\gamma C^1_\Phi h}
        \eqsp.
    \end{align}
\end{proposition}

\paragraph{Control problem.}

Let $\psi: \R^d\times \Nc(\R^d)\times A \to \R$ and $\Psi: \Nc(\R^d) \to \R$ be continuous functions, and consider the following assumption.
\begin{hypH}
\label{hypH:coercivity_hyp}
    Suppose that there exists $C_\Psi, c_\psi >0$ such that
    \begin{gather}
        - C_\Psi \left(1+\int_{\R^d}|y| \lambda(dy) + \langle 1\eqsp,\lambda\rangle\right)
        \leq
        \Psi(\lambda)
        \leq
        C_\Psi \left(1+\int_{\R^d}|y|^2 \lambda(dy) + \langle 1,\lambda\rangle^2 \right)\eqsp,
        \label{eq:coercivity_hyp:big_Psi}
        \\
        -C_\Psi \left(1+|x| \right)  + c_\psi|a|^2
        \leq
        \psi(x,\lambda, a )
        \leq
        C_\Psi \left(1+|x|^2 + 
        |a|^2 \right)
        \eqsp,
        \label{eq:coercivity_hyp:psi}
    \end{gather}
    for $(x,\lambda, a )\in \R^d\times \Nc(\R^d)\times A$.
\end{hypH}

Fix a standard strong control $\beta\in\CtrlStandard$ and a starting condition $(t,\check\lambda) \in[0,T]\times E$. The cost and  value functions are defined as follows:
\begin{align}
\label{eq:def:cost_function}
    J(t,\check\lambda;\beta) := \E^{\P}\left[
        \int_t^T \sum_{i\in \Vc^{t,\check\lambda;\beta}_s} \psi\left(Y^{i,\beta}_s,\xi^{t,\check\lambda;\beta}_s,\beta_s^i\right)ds +  \Psi\left(\xi^{t,\check\lambda;\beta}_T\right)\Bigg|\xi^{t,\check\lambda;\beta}_t=\check\lambda
    \right]
    \qquad
    \text{ and }
    \qquad
    v(t,\check\lambda) := \inf_{\beta\in\CtrlStandard} J(t,\check\lambda;\beta)
    \eqsp.
\end{align}

The well-posedness of this control problem follows directly from Proposition 2.4 in \citet{ocello2025controlled1}.

\paragraph{Symmetric controls.}

Consider now the following class of \emph{symmetric} controls, which assign same action to all particles in the same position.

\begin{definition}[Symmetric control]
    Fix $(t,\check\lambda) \in [0,T]\times E$.
    We say that $\beta = (\beta^i)_{i\in \Ic}$ is an \emph{symmetric control}, and we denote $\beta\in\SymmCtrl_{(t,\check\lambda)}$, if $\beta\in\CtrlStandard$ and, for $\xi^{t,\check\lambda;\beta}=\sum_{i\in\Vc^{t,\check\lambda;\beta}}\delta_{(i,Y^{i,\beta}_s)}$ solution of \eqref{SDE:strong}, we have
    \begin{align}
        \label{eq:cyclical_strong_ctrl_condition}
        \beta^i_s = \beta^j_s\eqsp,
        \quad \text{ whenever } Y^{i,\beta}_s = Y^{j,\beta}_s\eqsp,
    \end{align}
    for $s\in[0,T]$ and $i,j\in\Vc^{t,\check\lambda;\beta}_s$.
\end{definition}
 
As shown in \citet{ocello2025controlled1}, the invariance of the coefficients in the MF setting propagates to the feedback optimizer, making it sufficient to restrict the analysis to symmetric admissible controls. This reduction is formalized in the following proposition.

\begin{proposition}[Proposition 4.2 in \citet{ocello2025controlled1}]
    \label{prop:MF-symmetric-controls}
    Suppose Assumption H\ref{hypH:model_parameters}-H\ref{hypH:coercivity_hyp} hold. Fix $(t,\check\lambda=\sum_{i\in\Vc}\delta_{(i,x_i)}) \in [0,T]\times E$. Then, we have
    \begin{align}
        \label{eq:MF-symmetric-controls}
        v(t,\check\lambda) &= v(t,\mathfrak{s}\cdot\check\lambda)
        \eqsp,
        \quad\text{ for }\mathfrak{s}\in\mathfrak{S}_\Vc
        \eqsp,
        \qquad\text{ and }\qquad
        v(t,\check\lambda) = \inf_{\beta\in\SymmCtrl_{(t,\check\lambda)}} J(t,\check\lambda;\beta)
        \eqsp.
    \end{align}
\end{proposition}

For $t\in[0,T]$, $\lambda\in\Nc(\R^d)$, let $\valueStrong$ defined as 
\begin{align}
\label{eq:definition_v_symmetric}
    \valueStrong(t,\lambda) := 
    v\l(t,\check\lambda\r)
    \eqsp,
    \quad \text{ for any } \check\lambda\in E \text{ such that } \pi(\check\lambda) = \lambda
    \eqsp.
\end{align}
By Proposition \ref{prop:MF-symmetric-controls}, $\valueStrong$ is well-defined for any representation $\lambda$.
Therefore, we can restrict the state space of our analysis. Rather than working with the full configuration space $E$, we reformulate the problem directly on \(\Nc(\R^d)\), the space of empirical measures of particle positions. 
For $(t,\lambda)\in [0,T]\times\Nc(\R^d)$, with a slight abuse of notation, we continue to write $\xi^{t,\lambda;\beta}$ and say it satisfies \eqref{SDE:strong}, although what we truly mean is the projection $\pi(\xi^{t,\check\lambda;\beta}) \in \Nc(\R^d)$ of the original process $\xi^{t,\check\lambda;\beta}$ for some $\check\lambda\in E$ such that $\pi(\check\lambda) = \lambda$.

\paragraph{Pathwise uniqueness for the projected process.}

From \Cref{prop:MF-symmetric-controls}, our analysis can now be restricted to the projected process $\pi(\xi^{t,\hat\lambda;\beta})$, which captures the evolution of the empirical distribution of particle positions in the MF regime. The next step is to establish \emph{pathwise uniqueness} for this projected process. This property will play a crucial role in the sequel: it provides the foundation for connecting the relaxed formulation of the control problem with the strong one. In particular, it will allow us to show in \Cref{thm:equivalence_forms} that the two formulations are in fact equivalent.

To establish pathwise uniqueness, first consider the set $\mathfrak{S}^\text{tree}_\Vc\subset\mathfrak{S}_\Vc$ of permutations that preserve the genealogical structure of the branching process, originated from $\Vc$. More precisely, a permutation $\mathfrak{s}\in\mathfrak{S}^\text{tree}_\Vc$ if for any $i,j\in\Vc$ such that $k\preceq i,j$ for some $k\in\Vc$, we have $i\preceq j$ if and only if $\mathfrak{s}(i)\preceq \mathfrak{s}(j)$. Note that $\mathfrak{S}^\text{tree}_\Vc$ is a subgroup of $\mathfrak{S}_\Vc$.

\begin{corollary}
    \label{corollary:pathwise-uniqueness-projection}
    Suppose Assumption H\ref{hypH:model_parameters} holds. Fix $t\in [0,T]$, $\check\lambda= \sum_{i\in\Vc}\delta_{(i,x_i)} \in E$, and $\beta\in\SymmCtrl_{(t,\check\lambda)}$. Then, we have that
    \begin{align}
    \label{eq:pathwise-uniqueness-projection}
        \xi^{t,\check\lambda;\beta} = \mathfrak{s}^{-1}\cdot {}^{\mathfrak{s}}\xi^{t,\mathfrak{s}\cdot\check\lambda;\beta}
        \eqsp, \qquad\P-a.s., \quad \text{ for any }\mathfrak{s}\in\mathfrak{S}^\text{tree}_\Vc
        \eqsp,
    \end{align}
    where ${}^{\mathfrak{s}}\xi^{t,\mathfrak{s}\cdot\check\lambda;\beta}$ is the unique strong solution of \eqref{SDE:strong} with initial condition $(t, \mathfrak{s}\cdot\check\lambda)$ and control $\beta$, driven by the same family of noises $\{B^{\mathfrak{s}^{-1}(j)},Q^{\mathfrak{s}^{-1}(j)}\}_{j\in\Ic}$.
    Consequently, pathwise uniqueness holds for the projected process $\pi(\xi)$.
\end{corollary}
\begin{proof}
    The processes $\xi^{t,\check\lambda;\beta}$ and ${}^{\mathfrak{s}}\xi^{t,\mathfrak{s}\cdot\check\lambda;\beta}$ differ only in their initial labeling: one starts from $\check\lambda$, while the other starts from the relabeled configuration $\mathfrak{s}\cdot\check\lambda$.
    Nevertheless, by construction each particle branch is driven by the same source of randomness: for every \(i \in \Vc^{t,\check\lambda;\beta}_s\), the Brownian motion $B^i$ and Poisson random measure $Q^i$ that govern its evolution coincide with those driving the branch \(\mathfrak{s}(i) \in \Vc^{t,\mathfrak{s}\cdot\check\lambda;\beta}_s\).
    Thus, the two processes evolve under the same family of noises, consistently matched across relabeling of indices.

    Let $\{T_n\}_{n\ge0}$ be the (almost surely locally finite) sequence of jump times of $\xi^{t,\check\lambda;\beta}$ (and thus of$ {}^{\mathfrak{s}}\xi^{t,\mathfrak{s}\cdot\check\lambda;\beta}$, since both are driven by the same marked Poisson families and accept/reject the same marks under the previous alignment), with $T_0:=t$. Assume inductively that \eqref{eq:pathwise-uniqueness-projection} holds on $[t,T_n]$ for some $n\ge0$; we will show that it also holds on $[t,T_{n+1}]$.
    On the stochastic interval $[T_n,T_{n+1})$, no jump occurs, the coefficients are label-symmetric and depend only on the common projection, and each matched branch is driven by the same Brownian motion. Hence, pathwise uniqueness for the inter-jump SDEs implies  \eqref{eq:pathwise-uniqueness-projection} on $[T_n,T_{n+1})$. In particular their projections coincide on $[T_n,T_{n+1})$. At the jump time $T_{n+1}$, the thinning/marking construction uses the same Poisson marks for the matched labels, and since pre-jump states coincide, the acceptance decision and offspring configuration coincide as well; thus the post-jump projections are equal at $T_{n+1}$. By induction over $n$, we conclude that \eqref{eq:pathwise-uniqueness-projection} on the whole interval $[0,T]$. This argument is identical in spirit to the proof of Proposition 4.2 in \citet{ocello2025controlled1}.
    
    Finally, since $\pi$ is invariant under label permutations, we have $\pi(\xi^{t,\check\lambda;\beta}) = \pi({}^{\mathfrak{s}}\xi^{t,\mathfrak{s}\cdot\check\lambda;\beta})$, which establishes pathwise uniqueness for the projected process.
\end{proof}

\section{Relaxed formulation}
\label{section:relax_form}
We can give the \emph{relaxed formulation} for the branching diffusion control problem \eqref{eq:def:cost_function} by working with relaxed controls and weak solutions of the previous SDE \eqref{SDE:strong}, in the MF regime.

We equip the product space $[0,T] \times \R^d \times A$ with the $\sigma$-algebra $\Bc([0,T]) \otimes \Bc(\R^d) \otimes \Bc(A)$. Let $\Leb$ denote the Lebesgue measure in $\R_+$. Denote $\Ac^{\Leb}\subseteq M^1([0,T] \times \R^d \times A)$ the set of measures, whose projection on $[0,T]$ is absolutely continuous w.r.t.\ the Lebesgue measure $\Leb$.
Each $\alpha \in \Ac^{\Leb}$ can be identified with its disintegration \citep[see, \eg, Corollary 1.26, Chapter 1,][]{book:KALLENBERG-RM}. In particular, we have
\[
    \alpha(ds,dx,da) = ds \, \y_s(dx) \, \bar\alpha_s(x,da)
    \eqsp,
\]
for a process $\left(\y_s(dx)\right)_{s\geq0}$ (resp. $\left(\bar\alpha_s(x,da)\right)_{s\geq0,x\in\R^d}$) taking values in the set of functions from $[0,T]$ (resp. $[0,T]\times\R^d$) to $M^1(\R^d)$ (resp. $\Pc^1(A)$).

We endow on $\Ac^{\Leb}$ the \emph{stable topology} relative to the control space $A$. This topology is the weakest one that makes the mappings
$$
    \alpha \longmapsto \int_{[0,T]\times\R^d}\int_A \phi(s,x,a)\,\alpha(ds\,dx\,da)
$$
continuous, for all bounded measurable functions $\phi:[0,T]\times\R^d\times A\to\R$ that are continuous in the control variable $a \in A$.
From Proposition 2.10 in \citet{jacod2006type}, this topology is metrizable, and Theorem 2.8 therein provides a useful criterion for sequential compactness in this setting.

For $\x=\left(\x_s\right)_s\in\Db^{d}$ fixed, in theterminology of \citet{elkaroui1988existence}, we denote the \emph{space of generalized actions} $\Ac^{\Leb, \x}$ as
\begin{align*}
    \Ac^{\Leb, \x}:=\left\{\alpha\in \Ac^{\Leb} : \alpha(ds, dx, da) = ds \x_s(dx) \bar\alpha_s(x,da) \text{ a.e. }s\in [0,T]\right\}\eqsp,
\end{align*}
which is \textit{weakly* closed}.
We equip $\Db^{d} \times \Ac^{\Leb}$ with the product topology.

\subsection{Martingale model}
Let $\Lc$ be the generator defined on the cylindrical functions $F_\varphi = F(\langle \varphi, \cdot\rangle)$, for $F\in C^{2}_b(\R)$ and $\varphi\in C^2_{b}(\R^d)$, as
\begin{align*}
    &\Lc F_\varphi (x,\lambda,a)
    \\
    &= F'(\langle \varphi, \lambda\rangle) L \varphi(x,\lambda,a) + \frac{1}{2}
    F''(\langle \varphi, \lambda\rangle)
    \left|D\varphi(x)\sigma(x,\lambda,a)\right|^2
    + \gamma(x,\lambda,a)\bigg(
    \sum_{k\geq 0} F\big(\left\langle\varphi,\lambda\right\rangle + (k-1)\varphi(x)\big)p_k(x,\lambda,a)
    - F_\varphi\left(\lambda\right)\bigg) \eqsp.
\end{align*}
For simplicity, we write $F'_\varphi(\lambda)$, for $F'(\langle \varphi, \lambda\rangle)$ and $F''_\varphi(\lambda)$, for $F''(\langle \varphi, \lambda\rangle)$. 
Moreover, for $\mathfrak{F} = \{\mathscr{F}_s\}_{s\geq0}$ a filtration, we denote $\hat{\mathfrak{F}} =\{\hat{\mathscr{F}}_s\}_{s\geq0}$ the filtration such that $\hat{\mathscr{F}}_s := \Bc(\R^d)\otimes\mathscr{F}_s$, for $s\geq0$.

\begin{definition}[Relaxed control]
    \label{Def:Relaxed_control}
    Fix $(t,\lambda) \in [0,T]\times \Nc(\R^d)$. We say that $\Cc$ is a \emph{relaxed control}, and we denote $\Cc\in\RelCtrl_{(t,\lambda)}$, if
\begin{align*}
    \Cc = \left(
    \Omega, \mathscr{F}, \Pb,\mathfrak{F} =\left\{\mathscr{F}_s\right\}_{s\geq0}, X, \alpha
    \right)\eqsp,
\end{align*}
where
\begin{enumerate}[(i)]
    \item $\left(\Omega, \mathscr{F}, \Pb\right)$ is a probability space with complete right-continuous filtration $\mathfrak{F}$;
    \item $X=\left(X_s\right)_{s\geq0}$ is an $\mathfrak{F}$-progressively measurable process living in $\Db^{d}$ such that $\Pb(X_t=\lambda)=1$;
    \item $\bar\alpha:[0,T]\times \R^d\times \Omega \to \Pc^1(A)$ is a $\hat{\mathfrak{F}}$-predictable process associated with $\alpha\in\Ac^{\Leb}$ such that $\Pb(\alpha\in \Ac^{\Leb, X})=1$, \ie,
    \begin{align*}
        &\Pb\big(\alpha(ds, dx, da) = ds X_s(dx) \bar\alpha_s(x,da) \text{ a.e. }s\in [0,T]\big)=1,\\
        &\E^\Pb\left[\int_t^T \int_{\R^d\times A}|a|\bar\alpha_s(x,da)X_s(dx)ds\right]< \infty\eqsp;
    \end{align*}
    \item for $F_\varphi = F(\langle \varphi, \cdot\rangle)$, with $F\in C^{2}_b(\R)$ and $\varphi\in C^2_{b}(\R^d)$, the process $M^{F_\varphi}$ is a $(\Pb,\mathfrak{F})$-martingale, with
    \begin{align}
    \label{MartPb:diff-F_f}
        M^{F_\varphi}_s := F_\varphi(X_s) - \int_t^s \int_{\R^d\times A} \Lc F_\varphi (x,X_u,a)\bar\alpha_u(x,da)X_u(dx)du
        \eqsp, \qquad\text{ for } s\geq t
        \eqsp.
    \end{align}
\end{enumerate}
\end{definition}

For $\Cc\in\RelCtrl_{(t,\lambda)}$, we are only interested in the time interval $[t,T]$. Therefore, the process $X_s$ and the control $\alpha_s$ can be redefined for $s\in[0,t)$ as $X_s=\lambda$ and $\alpha_s=\delta_{a_0}$, for some $a_0\in A$.

\paragraph{Representation of the controlled martingale problem.}

There are several equivalent viewpoints for representing a measure-valued process, reflecting the fact that different formulations may be more convenient depending on the context. This is analogous to the case of real-valued processes, where one may work either with linear and quadratic variations or, alternatively, with Laplace functionals. In the present setting, the following lemma provides such an equivalent characterization of~\eqref{MartPb:diff-F_f}. It can be seen as the counterpart, in the controlled framework, of classical results such as Lemma 1.10 in \citet{etheridge2000introduction}, Theorem 1.3 in \citet{Roelly:convergence_measure_val_processes}, or Théorème 3.1 in \citet{Roelly:Rouault:Branchements_Spatiaux}. This result is proven considering the quadratic variation of a martingale \citep[see, \eg, Chapter I-4e,][]{jacod2013limit}.

\begin{lemma}\label{Lemma:Mart_charact:bd}
    Given $(t,\lambda) \in[0,T]\times \Nc(\R^d)$,
    let $\Cc = (
        \Omega, \mathscr{F}, \Pb,\mathfrak{F} =\{\mathscr{F}_s\}_{s\geq0}, X, \alpha
    )$ be such that conditions \textit{(i), (ii),} and \textit{(iii)} in the~\Cref{Def:Relaxed_control} are satisfied. The following are equivalent.

    \begin{enumerate}[(i)]
        \item We have $\Cc\in\RelCtrl_{(t,\lambda)}$.
        \item For any $\varphi\in C^2_{b}(\R^d)$ such that $\varphi>\eps$, for some $\eps>0$ and $\sup_{\R^d} \varphi \leq 1$, the process $M^{\exp_{\log\varphi}}$ is a $(\Pb,\mathfrak{F})$-martingale, with
        \begin{align}
        \label{MartPb:diff-exp}
            &M^{\exp_{\log\varphi}}_s :=
            e^{\langle \log\varphi,X_s\rangle} -
            \int_t^s \int_{\R^d\times A}  \left(
            \frac{L \varphi(x,X_u,a) + \gamma(x,X_u,a)( \Phi(\varphi(x),x,X_u,a) - \varphi(x))}{\varphi(x)} \right)\bar\alpha_u(x,da)X_u(dx)e^{\langle \log\varphi,X_u\rangle}du
            \eqsp,
        \end{align}
        for $s\geq t$.
        
        \item For any $\varphi\in C^2_b(\R^d)$, the process $\bar M^{\varphi}$ is a $(\Pb,\mathfrak{F})$-martingale, with 
        \begin{align}
        \label{MartPb:diff-var_finie}
        \begin{split}
            \bar M^{\varphi}_s = \langle\varphi,X_t\rangle -& \int_t^s \int_{\R^d\times A}   L \varphi(x,X_u,a)\bar\alpha_u(x,da)X_u(dx)du\\-&
            \int_t^s \int_{\R^d\times A}   \gamma(x,X_u,a)\left(\partial_s \Phi(1,x,X_u,a) - 1\right)\varphi(x)
            \bar\alpha_u(x,da)X_u(dx)du~, \quad \text{ for } s\in[t,T]\eqsp,
        \end{split}
        \end{align}
        having quadratic variation 
        \begin{align}
        \label{MartPb:diff-quadratic_var}
            \begin{split}
                \left[\bar M^{\varphi}\right]_s =&\int_t^s \int_{\R^d\times A}   \bigg( \text{Tr}\left(\sigma\sigma^\top(x,X_u,a)D\varphi D\varphi^\top(x)\right) \\
                &+\gamma(x,X_u,a)\left(
                \partial^2_{ss} \Phi(1,x,X_u,a) - \partial_s \Phi(1,x,X_u,a) + 1\right)\varphi^2(x)\bigg)
                \bar\alpha_u(x,da)X_u(dx)du\eqsp, \quad \text{ for }  s\in[t,T]\eqsp.
            \end{split}
        \end{align}
    \end{enumerate}
\end{lemma}

\begin{proof}

    $(i)\implies(ii)$: We need to prove that~\eqref{MartPb:diff-F_f} is a well defined martingale for the function $F_{\log\varphi}$ with $F(x)=\exp(x)$ and $\varphi\in C^2_{b}(\R^d)$ such that $\varphi>\eps$, for some $\eps>0$ and $\sup_{\R^d} \varphi \leq 1$. The process
    $M^{\exp_{\log\varphi}}$, as in~\eqref{MartPb:diff-exp}, is a local martingale. To prove that it is a martingale, we show its quadratic variation has a finite expectation. Since the compensator of $(M^{\exp_{\log\varphi}})^2$ is the same of $M^{\exp_{2\log\varphi}}=M^{\exp_{\log\varphi^2}}$, we get the quadratic variation of $M^{\exp_{\log\varphi}}$ applying~\eqref{MartPb:diff-F_f} to $F\in C^2_b(\R)$ and $\varphi^2$. Therefore, it is equal to
    \begin{align*}
    &\left[M^{\exp_{\log\varphi}}\right]_s =\int_t^s \int_{\R^d\times A}  \left(
        \frac{L \varphi^2(x,X_u,a) + \gamma(x,X_u,a)( \Phi(\varphi^2(x),x,X_u,a) - \varphi^2(x))}{\varphi^2(x)} \right)\\
        &\phantom{\left[M^{\exp_{\log\varphi}}\right]_s =\int_t^s \int_{\R^d\times A}L \varphi^2(x,X_u,a) + \gamma(x,X_u,a) \varphi^2(x)}
        \bar\alpha_u(x,da)X_u(dx) ~e^{\langle \log\varphi^2,X_u\rangle}du\eqsp.
    \end{align*}
    Since $\left[M^{\exp_{\log\varphi}}\right]$ is uniformly bounded, using Itô's isometry, $M^{\exp_{\log\varphi}}$ is a martingale.

    $(ii)\implies(iii)$:
    Fix $f\in C^2_{b}(\R^d)$. For $\theta>0$, and $M_f:=\sup_{\R^d}|f|$, we define $\varphi_1 := e^{\theta (f-M_f)}$ and $\varphi_2 := e^{-\theta M_f}$. Since $f$ is bounded, there exists $\eps>0$ such that $\varphi_1>\eps$ and $\sup_{\R^d}\varphi_1 \leq 1$. Applying~\eqref{MartPb:diff-exp} to $\varphi_1$ and $\varphi_2$, we get
    \begin{align}
        \label{eq:Mart_charact:bd:varphi1}
        \begin{split}
            &\E^{\Pb}\bigg[e^{\langle \theta (f-M_f),X_{s+h}\rangle} - e^{\langle \theta (f-M_f),X_{s}\rangle}
            \\
            &\quad- \int_{s}^{s+h} \int_{\R^d\times A} \bigg(\theta L f(x,X_u,a) +
            \theta^2 \text{Tr}\left(\sigma\sigma^\top(x,X_u,a)Df Df^\top(x)\right)
            \\
            &\quad\qquad\qquad\qquad\qquad
            + \gamma(x,X_u,a) \frac{ 
                \Phi\left(\left(e^{\theta (f(x)-M_f)}\right),x,X_u,a\right) - e^{\theta (f(x)-M_f)}
            }{e^{\theta (f(x)-M_f)}}\bigg)
            \bar \alpha_u(x,da) X_u(dx) e^{\langle \theta (f-M_f),X_{u}\rangle}du\bigg| \Fc_s\bigg]=0 
            \eqsp,
        \end{split}
        \\
        \label{eq:Mart_charact:bd:varphi2}
        \begin{split}
            &\E^{\Pb}\bigg[e^{\langle -\theta M_f,X_{s+h}\rangle} - e^{\langle -\theta M_f,X_{s}\rangle}
            \\
            &\quad- \int_{s}^{s+h} \int_{\R^d\times A} \gamma(x,X_u,a)
            \frac{\Phi\left(e^{-\theta M_f},x,X_u,a\right) - \left(e^{-\theta M_f}\right)}{e^{-\theta M_f}}
            \bar \alpha_u(x,da)X_u(dx)e^{\langle -\theta M_f,X_{u}\rangle}du\bigg| \Fc_s\bigg]=0 
            \eqsp.
        \end{split}
    \end{align}
    Since all the functions are bounded, we are allowed to differentiate with respect to $\theta$. Dividing by $\theta$, subtracting~\eqref{eq:Mart_charact:bd:varphi1} and~\eqref{eq:Mart_charact:bd:varphi2}, and setting $\theta=0$, we get~\eqref{MartPb:diff-var_finie}. Differentiating twice with respect to $\theta$, dividing by $\theta^2$subtracting~\eqref{eq:Mart_charact:bd:varphi1} and~\eqref{eq:Mart_charact:bd:varphi2} and setting $\theta=0$, we get~\eqref{MartPb:diff-quadratic_var}.

    $(iii)\implies(i)$: We prove the last implication using Itô's formula for semimartingales. Fix $F\in C^2(\R^n)$ and $f\in C^2_b(\R^n)$. We have that $\langle f,X_{s}\rangle_{s\geq t}$ is a $\Pb$-semimartingale, and so, by Itô's formula, we have~\eqref{MartPb:diff-F_f}.
\end{proof}



\paragraph{Relaxed control problem.}
We can now define the relaxed control problem.
For $\Cc\in \RelCtrl_{(t,\lambda)}$, we define the cost function as
\begin{align}
\label{control_pb:relaxed:cost_fct}
    J(t,\lambda;\Cc) = 
    \E^{\Pb}\left[
    \int_t^T \int_{\R^d\times A} \psi\left(s,X_s,a\right)\bar\alpha_s(x,da)X_s(dx)ds +  \Psi\left(X_T\right)
    \right]\eqsp,
\end{align} 
and the \textit{relaxed control problem} as
\begin{align}\label{control_pb:relaxed:value_fct}
    \valueRelaxed(t,\lambda) := \inf\left\{J(t,\lambda;\Cc)~:~\Cc\in \RelCtrl_{(t,\lambda)}\right\}\eqsp,
    \qquad\text{ for }(t,\lambda) \in [0,T]\times \Nc(\R^d)
    \eqsp.
\end{align}

Note that the \emph{relaxed formulation} of the control problem is \emph{well-posed}. In particular, using the martingale problem \eqref{MartPb:diff-F_f} in place of the strong Itô's formula, we can deduce finiteness of moments from the coercivity properties of the cost functional in Assumption H\ref{hypH:coercivity_hyp}. This follows from arguments that are direct generalizations of Proposition 2.2, Lemma 2.3, and Proposition 2.4 in \citet{ocello2025controlled1}. Therefore, the extension of these results is immediate and their proofs are omitted.

\begin{remark}
\label{rmk:embedding-symmetric-controls}
    It is straightforward to embed the class of strong controls into the relaxed framework. Specifically, for any $\check\lambda \in E$ with $\pi(\check\lambda)=\lambda$, we have the inclusion
    \begin{align*}
        \CtrlStandard^\mathfrak{s}{(t,\check\lambda)} \subseteq \Rc^\mathfrak{r}{(t,\lambda)}
        \eqsp.
    \end{align*}
    Indeed, given $\beta \in \CtrlStandard^\mathfrak{s}_{(t,\check\lambda)}$ with associated controlled process $\xi^{t,\check\lambda;\beta}$, we define the relaxed state-control pair by
    \begin{align*}
        X_s := \pi\bigl(\xi^{t,\check\lambda;\beta}_s\bigr)
        \qquad
        \text{ and }
        \qquad
        \alpha(ds,\,dx,\,da) := ds\,\pi\bigl(\xi^{t,\check\lambda;\beta}_s\bigr)(dx)\,\delta_{\beta_s^i(x)}(da),
    \end{align*}
    for $s\in[0,T]$,
    where $\delta{\beta^i_s(x)}$ denotes the Dirac mass at the control $\beta^i_s \in A$ corresponding to the unique particle located at $x$, a point in the support of $\pi\bigl(\xi^{t,\check\lambda;\beta}_s\bigr)$. This construction is well defined due to the symmetry of the control. Consequently, the infimum in the relaxed formulation \eqref{control_pb:relaxed:value_fct} cannot exceed that of the strong formulation \eqref{eq:def:cost_function}.
\end{remark}

We aim at proving the converse inclusion, thus the equivalence between the strong and the relaxed formulations of the control problem, \ie, $\valueRelaxed(t,\lambda) = v(t,\lambda)$, for $(t,\lambda)\in[0,T]\times\Nc(\R^d)$. To do so, we follow the path set by \citet{elkaroui1987compactification} and \citet{haussmann1990existence}, adapting their results to our branching framework.


\subsection{Natural controls}
\label{Subsection:natural-ctrl}

We state the following straightforward adaptation of Lemma 3.7 of \citet{haussmann1990existence}. This construction allows the process $X$ to be considered with respect to its canonical filtration. We stress that the following lemma is stated relative to the raw filtration generated by the processes, without taking the right-continuous modification or the completion under a given probability measure. This choice is consistent with the framework of \citet{haussmann1990existence}, where the key requirement is simply the existence of a countable, dense family of test functions that characterizes the associated martingale problem.

\begin{lemma}\label{Lemma:F_X-measurability}
Fix $(t,\lambda) \in [0,T]\times \Nc(\R^d)$ and $\Cc= (
\Omega, \mathscr{F}, \Pb, \{\mathscr{F}_s\}_{s\geq 0}, X, (\bar\alpha_s)_{s\geq 0}
)\in\RelCtrl_{(t,\lambda)}$. If $\{\mathscr{F}^X_s\}_{s\geq0}$ is the filtration generated by $X$ and $\{\mathscr{G}_s\}_{s\geq 0}$ another filtration such that $\mathscr{F}^X_s\subseteq\mathscr{G}_s\subseteq\mathscr{F}_s$, for $s\geq 0$. Then, there exists $\alpha^\mathscr{G}$ such that 
\begin{align*}
\bar\Cc = \left(
\Omega, \mathscr{G}_T, \Pb, \left\{\mathscr{G}_s\right\}_{s\geq 0}, X, \alpha^\mathscr{G}
\right)
\end{align*}
is in $\RelCtrl_{(t,\lambda)}$ and $J(t,\lambda;\Cc)=J(t,\lambda;\bar\Cc)$.
\end{lemma}

We now introduce the class of \emph{natural controls}, following the terminology of \citet{krylov1980controlled} and later adopted in \citet{haussmann1990existence}.
Let $\mu$ be the canonical process on $\Db^{d}$, and denote by $\F^{\mu}=\left\{\Fc_s^{\mu}\right\}_{s\geq 0}$ its right continuous filtration.
Motivated by the previous lemma, we restrict attention to the following subclass of relaxed controls.

\begin{definition}
Fix $(t,\lambda) \in[0,T]\times \Nc(\R^d)$. $\Cc= (
\Omega, \mathscr{F}, \Pb, \left\{\mathscr{F}_s\right\}_s, X, \alpha)$ in $\RelCtrl_{(t,\lambda)}$ is a \emph{natural control}, and we say that $\Cc$ is in $\NatCtrl_{(t,\lambda)}$, if $\Omega=\Db^{d}$, $\mathscr{F} = \Fc^\mu_T$, $\mathscr{F}_s = \Fc^\mu_s$, for $s\in[t,T]$, $X=\mu$, and
\begin{align*}
    \Pb\left(\mu_s=\lambda, s\in [0,t]\right)=1\eqsp.
\end{align*}
\end{definition}

We note that, in the previous definition, the only objects not predetermined are the pair $(\Pb,\alpha)$, consisting of a probability measure $\Pb$ on $\Db^d$ (the law of the canonical process $\mu$) and the control process $\alpha$. With a slight abuse of notation, we identify $(\Pb,\alpha)$ with the natural control system
$\Cc^{\Pb,\alpha}:=(
\Db^{d}, \Fc^\mu_T, \Pb, \{\Fc^\mu_s\}_{s\geq0}, \mu, \alpha)$ in $\NatCtrl_{(t,\lambda)}$.

\paragraph{Weak controls.}
As noted in \Cref{rmk:embedding-symmetric-controls}, the class of strong controls in embedded in the relaxed controls and this embedding corresponds to the class of relaxed controls where the measure-valued process is almost surely a Dirac measure. It is therefore natural to consider the subclass of natural controls where the measure-valued process is almost surely a Dirac measure. This class is known in the literature as \emph{weak controls} \citep[see, \eg,][]{haussmann1990existence}.


For a fixed $\x\in\Db^{d}$, the set of measurable functions $\mathfrak{a}:[0,T]\times \R^d\to A$ is canonically embedded in $\Ac^{\Leb, \x}$ by $\alpha^\mathfrak{a}(ds,dx, da):= ds \eqsp \x_s(dx)\eqsp \delta_{\mathfrak{a}(s,x)}(da)$.

\begin{definition}
    Fix $(t,\lambda) \in[0,T]\times \Nc(\R^d)$. We say that $(\Pb,\mathfrak{a})$ is an \emph{weak control}, and we write $(\Pb,\mathfrak{a})\in\WeakCtrl_{(t, \lambda)}$, if $\mathfrak{a}:[0,T]\times \R^d\times \Db^d \to A$ is $\hat\F^\mu$-predictable, and $\left(\Pb, \alpha^\mathfrak{a}\right) \in \NatCtrl_{(t, \lambda)}$.
\end{definition}

Therefore, for $\Pb\in\WeakCtrl_{(t, \lambda)}$, we have that
\begin{align*}
    F_\varphi(\mu_s) - \int_t^s \int_{\R^d} \Lc F_\varphi (x,\mathfrak{a}(u,x),\mu_u)\mu_u(dx)du
\end{align*}
is a $(\Pb, \F^{\mu})$-martingale for $s\geq t$, $F\in C^{2}_b(\R)$ and $\varphi\in C^2_{b}(\R^d)$.

\subsection{Filippov's convexity condition}
We now introduce a condition under which the class of admissible controls can be restricted from \(\NatCtrl_{(t,\lambda)}\) to \(\WeakCtrl_{(t,\lambda)}\) without altering the value function, for $\varepsilon$-optimal controls. This relies on showing that, under the following assumption, one can always associate to each natural control a weak control yielding the same cost.

\begin{hypH}\label{Assumption:Convex}
The following set
\begin{align*}
    K(x,\lambda):=\left\{\left(b(x,\lambda,a),\sigma\sigma^\top(x,\lambda,a),\big((\gamma p_k)(x,\lambda,a)\big)_{k\geq0}, z\right) : a\in A, z\geq \psi(x,\lambda,a) \right\}
    \subseteq \R^d\times\R^{d\times d}\times \R_+^\infty\times \R
\end{align*}
is convex for all $(x,\lambda)\in\R^d\times \Nc(\R^d)$.
\end{hypH}

This convexity assumption is the so-called \emph{Filippov condition}, which is common in the control literature. It holds when $A$ is a convex subset of a vector space, and the parameters are affine in $a$, which is the case of the linear-quadratic examples in \citet{ocello2025controlled1}.

\begin{proposition}
\label{Prop:atom_ctrlVSnat_ctrl}
    Fix $(t,\lambda) \in[0,T]\times \Nc(\R^d)$ and $\varepsilon>0$. Suppose that Assumptions H\ref{hypH:model_parameters}--H\ref{hypH:coercivity_hyp}--H\ref{Assumption:Convex} hold. For $(\Pb, \alpha) \in \NatCtrl_{(t, \lambda)}$ such that
    \begin{align*}
        J\l(t,\lambda;\Cc^{\Pb,\alpha}\r)
        \leq
        \inf\left\{
            J\l(t,\lambda;\Cc^{\tilde{\Pb},\tilde\alpha}\r) ~:~(\tilde\Pb, \tilde\alpha) \in \NatCtrl_{(t, \lambda)}
        \right\} + \varepsilon \eqsp,
    \end{align*}
    there exists $\mathfrak{a}$ such that $(\Pb, \mathfrak{a}) \in \WeakCtrl_{(t, \lambda)}$ and $J(t,\lambda;\Cc^{\Pb,\alpha^\mathfrak{a}})\leq J(t,\lambda;\Cc^{\Pb,\alpha})$.
\end{proposition}

\begin{proof}
    This proof is an adaptation of the proof of Theorem 3.6 of \citet{haussmann1990existence}.

    Fix $(\Pb, \left(\alpha_s\right)_s)\in\NatCtrl_{(t, \lambda)}$. Define $c^1$ and $c^2$ as
    \begin{align*}
        c^1(s,x,\omega) := \int_{A} \left(b, \sigma\sigma^\top,(\gamma p_k)_{k\geq0}\right)(x,\mu_s(\omega),a)\bar \alpha_s(x, da)(\omega)
        \eqsp,\qquad
        c^2(s,x,\omega) := \int_{A} \psi(x,\mu_s(\omega),a)\bar \alpha_s(x, da)(\omega)
        \eqsp.
    \end{align*}
    Since all the functions defining $K$ are continuous, it follows that for almost every $(x,\lambda)\in\R^d\times \Nc(\R^d)$, the set $K(x,\lambda)$ is closed. This result can be seen as a direct adaptation of Proposition 3.5 in \citet{haussmann1990existence}.
    Combining this with Assumption H\ref{Assumption:Convex}, we have that $K(x,\lambda)$ is closed and convex. Therefore, $(c^1,c^2)(s,x,\omega)$ is in $K(x,\mu_s(\omega))$, for $\Leb\otimes\mu_s(\omega)$-almost all $(s, x)$, and $\Pb$-almost all $\omega$. Moreover, $(c^1,c^2)$ is $\hat\F^\mu$-predictable.
    
    Applying Theorem A.9 of \citet{haussmann1990existence}, there is a $\hat\F^\mu$-predictable $A$-valued process $\mathfrak{a}$ such that
    \begin{align}
        \label{eq:prop_atomVSnat:meas_sel1}
        c^1(s,x,\omega) = \left(b, \sigma\sigma^\top,(\gamma p_k)_{k\geq0}\right)(x,\mu_s(\omega),\mathfrak{a}(s,x,\omega))
        \eqsp,
        \qquad
        c^2(s,x,\omega)\geq \psi(x,\mu_s(\omega),\mathfrak{a}(s,x,\omega))
        \eqsp,
    \end{align}
    for $\Leb\otimes\mu_s(\omega)$-almost all $(s, x)$, and $\Pb$-almost all $\omega$.
    For $F\in C^{2}_b(\R)$ and $\varphi\in C^2_{b}(\R^d)$, we must have
    \begin{align*}
        \int_{\R^d\times A} \Lc F_\varphi (x,\mu_u(\omega),a)\bar\alpha_u(x,da;\omega)\mu_u(dx;\omega) = \int_{\R^d} \Lc F_\varphi (x,\mu_u(\omega),\mathfrak{a}(s,x,\omega))\mu_u(dx;\omega)\eqsp,
    \end{align*}
    for $\Leb\otimes\mu_s(\omega)$-almost all $(s, x)$, and $\Pb$-almost all $\omega$. This means that the process
    $$
        F_\varphi(\mu_s) - \int_t^s \int_{\R^d\times A} \Lc F_\varphi (x,\mathfrak{a}(s,x,\mu_u),\mu_u)\mu_u(dx)du
    $$ is a $\Pb$-martingale, for $s\geq t$. This means that the law of $(\mu_s)_{s\geq 0}$ under the control $(\Pb,{\alpha})$ is the same as the law of the solution of the SDE associated to the weak control $(\Pb,\mathfrak{a})$.

    Combining Assumption H\ref{hypH:coercivity_hyp} with \eqref{eq:prop_atomVSnat:meas_sel1}, we have that
    \begin{align*}
        &-C_\Psi\eqsp
        \E^\Pb\left[ 
        \int_t^T \left(
            \langle1,\mu_s\rangle + \int_{\R^d} |x|\mu_s(dx)
        \right)ds
        \right]
        + c_\psi\eqsp
        \E^\Pb\left[
            \int_t^T \int_{\R^d} |\mathfrak{a}(s,x,\cdot)|^2\mu_s(dx)ds
        \right]
        \\
        &\leq
        \E^\Pb\left[
            \int_t^T \int_{\R^d\times A} \psi(x,\mu_s,\mathfrak{a}(s,x,\cdot))\mu_s(dx)ds
        \right]
        \\
        &
        \leq
        \E^\Pb\left[
            \int_t^T \int_{\R^d\times A} \psi(x,\mu_s,a)\bar\alpha_s(x,da)\mu_s(dx)ds
        \right]
        \\
        &\leq
        J(t,\lambda;\Cc^{\Pb,\alpha}) - 
        \E^\Pb\left[
            \Psi(\mu_T)
        \right]
        \\
        &\leq
        \inf\left\{
            J\l(t,\lambda;\Cc^{\tilde{\Pb},\tilde\alpha}\r) ~:~(\tilde\Pb, \tilde\alpha) \in \NatCtrl_{(t, \lambda)}
        \right\}
        +\varepsilon - 
        \E^\Pb\left[
            \Psi(\mu_T)
        \right]
        \eqsp.
    \end{align*}
    Using Assumption H\ref{hypH:coercivity_hyp} again, there exists $C>0$ such that
    \begin{align*}
        \E^\Pb\left[
            \int_t^T \int_{\R^d} |\mathfrak{a}(s,x,\cdot)|^2\mu_s(dx)ds
        \right]
        \leq C \eqsp T \eqsp  
        \E^\Pb\left[
            \sup_{s\in[t,T]} \langle1, \mu_s\rangle^2+ \langle|\cdot|^2, \mu_s\rangle
        \right]+
        \inf\left\{
            J\l(t,\lambda;\Cc^{\tilde{\Pb},\tilde\alpha}\r) ~:~(\tilde\Pb, \tilde\alpha) \in \NatCtrl_{(t, \lambda)}
        \right\}
        +\varepsilon
        \eqsp.
    \end{align*}
    Following the steps of the proof of Proposition 
    2.2 and Proposition 2.4 in \citet{ocello2025controlled1}, since $(\Pb,\alpha)$ is an $\varepsilon$-optimal control, we have that we can bound uniformly the right-hand side of the previous inequality. This means that, in particular, we have
    \[
        \E^\Pb\left[\int_t^T \int_{\R^d\times A}|\mathfrak{a}(s,x,\cdot)|\mu_s(dx)ds\right]< \infty
        \eqsp,
    \]
    which entails that $(\Pb,\mathfrak{a})\in\WeakCtrl_{(t,\lambda)}$. Moreover, from~\eqref{eq:prop_atomVSnat:meas_sel1}, we get $J(t,\lambda;\Cc^{\Pb,\alpha^\mathfrak{a}})\leq J(t,\lambda;\Cc^{\Pb,\alpha})$.
\end{proof}

\section{Existence of an Optimal Control}
\label{section:Existence_optimum}

\subsection{Equivalence between Relaxed and Strong Formulations}
\label{Section:equivalence-formulations}
We are now in a position to establish the equivalence between the strong and relaxed formulations of the control problem.

\begin{theorem}
\label{thm:equivalence_forms}
Suppose Assumption H\ref{hypH:model_parameters}--H\ref{hypH:coercivity_hyp}--H\ref{Assumption:Convex} hold. For $(t,\lambda)\in[0,T]\times \Nc(\R^d)$, we have
\begin{align}
\label{eq:equivalence_forms}
    \valueStrong(t,\lambda) = \inf\left\{J(t,\lambda;\Cc):\Cc\in \RelCtrl_{(t,\lambda)}\right\}
    =
    \inf\left\{J(t,\lambda;\mathfrak{a}):\mathfrak{a}\in \WeakCtrl_{(t,\lambda)}\right\}
    \eqsp,
\end{align}
with $v$ defined in~\eqref{eq:definition_v_symmetric}.
\end{theorem}

\begin{proof}
    \emph{Step 1.}
    From \Cref{rmk:embedding-symmetric-controls}, we know that every strong control can be naturally embedded into the relaxed framework, so that the relaxed formulation always provides a lower bound on the strong value function, \ie,
    \begin{align*}
        \valueStrong(t,\lambda) \geq \inf\left\{J(t,\lambda;\Cc):\Cc\in \RelCtrl_{(t,\lambda)}\right\}
        \eqsp.
    \end{align*}

    \emph{Step 2.} We know show that we can focus on the control problem w.r.t.\ the set of weak controls. To this end, we adapt to our setting Corollary 3.9 in \citet{haussmann1990existence}.
    
    Fix $(t,\lambda)\in[0,T]\times \Nc(\R^d)$. Consider $\Cc\in\RelCtrl_{(t,\lambda)}$ such that
    \begin{align*}
        J(t,\lambda;\Cc) \leq \valueRelaxed(t,\lambda) + \varepsilon
        \eqsp,
    \end{align*}
    with
    \begin{align*}
        \Cc := \l(
            \Omega, \mathscr{F}, \Pb, \{\mathscr{F}_s\}_{s\geq 0}, X, (\bar\alpha_s)_{s\geq 0}
        \r)
        \eqsp.
    \end{align*}
    Applying \Cref{Lemma:F_X-measurability}, we can find a control $\tilde\Cc\in\RelCtrl_{(t,\lambda)}$ such that $J(t,\lambda;\tilde\Cc) = J(t,\lambda;\Cc)$ and
    \begin{align*}
        \tilde\Cc := \l(
            \Omega, \mathscr{F}, \Pb, \{\mathscr{F}^X_s\}_{s\geq 0}, X, \tilde\alpha
        \r)
        \eqsp.
    \end{align*}
    Moreover, from a straightforward generalization of \Cref{Prop:atom_ctrlVSnat_ctrl}, we can find a $\{\hat{\mathscr{F}}^X_s\}_{s\geq 0}$-predictable control process $\mathfrak{a}^X$ such that 
    \begin{align*}
        \tilde{\tilde{\Cc}} := \l(
            \Omega, \mathscr{F}, \Pb, \{\mathscr{F}^X_s\}_{s\geq 0}, X, ds X_s(dx)\delta_{\mathfrak{a}^X(s,x, \cdot)}
        \r)
    \end{align*}
    and $J(t,\lambda;\tilde{\tilde{\Cc}})\leq J(t,\lambda;\tilde{\Cc})$.
    Let now $\Pb^X$ be the law of $X$ as a process in $\Db^d$.
    From Doob's functional representation theorem \citep[see, \eg, Lemma 1.13,][]{kallenberg1997foundation}, there exists a $\Bc([0,T])\otimes\Bc(\R^d)\otimes \Bc(\Db^d)$-measurable function $\kappa^\mathfrak{a}: [0,T]\times \R^d \times \Db^d \rightarrow A$ such that
    \begin{align*}
        \mathfrak{a}^X(s,x,\omega)  =  \kappa^\mathfrak{a}(
            s,x,X_{s-\wedge\cdot}(\omega)
        )
        =
        \kappa^\mathfrak{a}(s,x,X(\omega))
        \eqsp, \qquad s\in[0,T],\eqsp x\in\R^d, \eqsp \omega\in \Omega\eqsp.
    \end{align*}
    Therefore, from Lemma A.1 in \citet{haussmann1990existence}, consequence of Lemma 7, Appendix 1 in \citet{Dellacherie:Meyer:B}, there exists a weak control $\mathfrak{a}^\mu$ such that $J(t,\lambda;\tilde{\tilde{\Cc}}) = J(t,\lambda;\Cc^{\Pb^X,\alpha^{\mathfrak{a}^\mu}})$.

    \emph{Step 3.}
    The last step consists in upgrading the weak control to a strong one. By the \emph{Yamada--Watanabe principle} from \citet{yamada1971uniqueness} \citep[see also Theorem~32.14 in][]{kallenberg1997foundation}, the combination of \emph{weak existence} together with \emph{pathwise uniqueness} implies the \emph{existence of a unique strong solution} adapted to the given noise. 

    This argument is not restricted to purely diffusive dynamics, as it has been extended to jump processes driven by Poisson random measures, where an analogous equivalence between weak existence plus pathwise uniqueness and strong existence holds \citep[see, \eg, Chapter 4 in][]{situ2005theory}. Moreover, our setting can be recovered as a particular case of the abstract framework developed in \citet{kurtz2007yamada}, which generalizes the Yamada–Watanabe principle to a broad class of stochastic models.

    Fix some $a_0\in A$. We consider the filtered space $(\Omega,\Fc,\F,\P)$ as in~\Cref{Section:strong_form}. We can then define the standard strong control $\beta^\mathfrak{a}$ as
    \begin{align*}
        \beta^{\mathfrak{a},i}_s
        = 
        \kappa^\mathfrak{a}\left(
            s,Y^{i,\beta}_{s-}, 
            \left(
                \xi^{t,\lambda;\beta}_{u\wedge s-}
            \right)_{u\in[0,T]}
        \right)
        \1_{i\in \Vc^{t,\lambda;\beta}_{s-}}
        + a_0\eqsp
        \1_{i\notin \Vc^{t,\lambda;\beta}_{s-}}\eqsp, \quad \text{ for } i\in \Vc^{t,\lambda;\beta}_{s-}\eqsp,
    \end{align*}
    where $\xi^{t,\lambda;\beta}$ (resp. $Y^{i,\beta}$, for $i\in \Vc^{t,\lambda;\beta}_s$) is the strongly controlled population (resp. particle) associated with $\beta^{\mathfrak{a}}$.
    From the construction of the control, we immediately obtain weak existence of the process $\xi$, whose law is by definition $\Pb^X$. Furthermore, by \Cref{corollary:pathwise-uniqueness-projection}, pathwise uniqueness holds for the projected dynamics in $\Db^d$ governed by \eqref{SDE:strong}. By the Yamada–Watanabe principle, weak existence together with pathwise uniqueness implies strong existence. Consequently, the strong control problem coincides with the formulation in terms of weak controls, thereby establishing \eqref{eq:equivalence_forms}.

\end{proof}
\subsection{Control rules}
\label{Section:ctrl-rules}

In order to finally prove the existence of an optimal control, it is sufficient, by \Cref{thm:equivalence_forms} together with \Cref{Lemma:F_X-measurability}, to work within any intermediate formulation $\Rc^{\cdot}$ satisfying
\begin{align*}
    \WeakCtrl \subseteq \Rc^{\cdot} \subseteq \RelCtrl\eqsp,
\end{align*}
since restricting to any of these classes does not alter the value function. This flexibility allows us to reformulate the problem directly in terms of \emph{canonical relaxed controls}, as done in \citet{haussmann1990existence} and \citet{lacker2017limit}, which are naturally defined on the product space $\Omega=\Db^d\times \Ac^{\Leb}$ with canonical processes $(\mu,\alpha)$ and filtration $\F^{\mu,\alpha}={\Fc_s^{\mu,\alpha}}_{s\geq 0}$ generated by them, \ie, the following filtration 
\begin{align*}
    \sigma\left(\mu_s(\Ec_1), \a([0,s']\times \Ec_2\times \Ec_3],\text{ for }s,s'\in[0,T],\eqsp \Ec_1,\Ec_2\in\Bc(\R^d),\eqsp\Ec_3\in\Bc(A)\right)\eqsp.
\end{align*}

\begin{definition}[Control rule]
Fix $(t,\lambda) \in[0,T]\times \Nc(\R^d)$.
$\Cc= (
\Omega, \mathscr{F}, \Pb, \{\mathscr{F}_s\}_{s\geq0}, X, \alpha)\in \RelCtrl_{(t,\lambda)}$ is a \emph{control rule}, and we write $\Cc\in\CtrlRule_{(t,\lambda)}$, if $\Omega=\Db^{d}\times \Ac^{\Leb}$, $\mathscr{F} = \Fc^{\mu,\a}_T$, $\mathscr{F}_s = \Fc^{\mu,\a}_s$, for $s\in[t,T]$, $X=\mu$, $\alpha=\a$, and
\begin{align*}
    \Pb\left(\mu_s=\lambda,\eqsp s\in [0,t]\right)
    =1
    \eqsp.
\end{align*}
\end{definition}

With abuse of notation, we write $\Pb\in\Rc_{(t,\lambda)}$ (resp. $J(t,\lambda;\Pb)$) to denote $\Cc^\Pb:=(
    \Db^{d}\times \Ac^{\Leb}, \Fc^{\mu,\a}_T, \Pb, \left\{\Fc^{\mu,\a}_s\right\}_{s\geq0}, \mu, \a
)\in\Rc_{(t,\lambda)}$ (resp. $J(t,\lambda;\Cc^\Pb)$). From~\Cref{Lemma:F_X-measurability}, any relaxed control is associated with a control rule with the same cost function $J$. Therefore, from~\Cref{thm:equivalence_forms}, we have
\begin{align*}
    \valueStrong(t,\lambda) =
    \valueRelaxed(t,\lambda) =
    \inf\left\{
        J(t,\lambda;\Pb)\eqsp:\eqsp\Pb\in \CtrlRule_{(t,\lambda)}
    \right\}
    \eqsp.
\end{align*}

This definition is particularly useful since a control rule is completely characterized by a probability measure $\Pb \in \Pc^1(\Db^d\times \Ac^{\Leb})$, which prescribes the joint distribution of the canonical processes $(\mu,\mathbf{a})$.
Consequently, the analysis of optimal controls reduces to studying the topological properties of the space $\Pc^1(\Db^d\times \Ac^{\Leb})$. In particular, our strategy is to show that the optimization problem amounts to minimizing a lower semicontinuous functional over a compact set. To this end, it suffices to establish that the cost functional $J$ is \emph{lower semicontinuous} and that, for each $\varepsilon>0$, the set $\CtrlRuleEps_{(t,\lambda)}$ is \emph{compact} in $\Pc^1(\Db^{d} \times \Ac^{\Leb})$, with 
\begin{align*}
    \CtrlRuleEps_{(t,\lambda)}
    := \left\{
        \Pb\in \CtrlRule_{(t,\lambda)} \eqsp:\eqsp J(t,\lambda;\Pb) \leq \valueRelaxed(t,\lambda)+ \eps
    \right\}
    \eqsp.
\end{align*} 

\subsection{Existence of an Optimal Control}

\paragraph{Lower semicontinuity.}
First, we establish the lower semi-continuity of the cost function.

\begin{lemma}\label{Lemma:J_continuous}
    For $(t,\lambda) \in[0,T]\times \Nc(\R^d)$, $J(t,\lambda;\cdot)$ is lower semicontinuous on $\Pc^1(\Db^{d} \times \Ac^{\Leb})$.
\end{lemma}
\begin{proof}
    Consider $f:\Db^{d} \times \Ac^{\Leb} \to \R$, defined as
    \begin{align*}
        f(\x,\alpha) := \int_t^T \int_{\R^d\times A}\psi\left(x,\x_s,a\right)\bar \alpha_s(x,da)\x_s(dx) ds +\Psi\left(\x_T\right)
        \eqsp.
    \end{align*}
    Define, for $(\x,\alpha)\in \Db^{d} \times \Ac^{\Leb}$, the measure $\Gamma^{x,\alpha}(ds,dy,da):=ds\,x_s(dy)\,\bar\alpha_s(y,da)$ on $[t,T]\times\R^d\times A$.
    Since $\psi$ is continuous, the map $(z:= (s,y,a))\mapsto \phi_{\x}(z):=\psi(y,\x_s,a)$ is lower semicontinuous on $[t,T]\times\R^d\times A$ and bounded from below up to the coercive growth \eqref{eq:coercivity_hyp:psi}. By the lower semicontinuity of integrals w.r.t.\ narrowly converging measures with l.s.c.\ integrand (Portmanteau's theorem), we get
    \begin{align*}
        \liminf_{n\to\infty}\int \psi\bigl(y,\x^n_s,a\bigr)\, \Gamma^{\x^n,\alpha^n}(ds,dy,da)\;\ge\;\int \psi\bigl(y,\x_s,a\bigr)\,\Gamma^{\x,\alpha}(ds,dy,da)
        \eqsp,
    \end{align*}
    for a sequence $\{(\x_n,\alpha_n)\}_{n\geq 0}$ 
    converging in $\Db^d\times \Ac^{\Leb}$ to $(\x,\alpha)$.
    For the terminal part, $\x\mapsto \x_T$ is continuous as a map $\Db^{d}\to\R^d$ at every $\x$ (and in particular at the limit point), hence $\Psi(\x^n_T)\to \Psi(\x_T)$ by continuity of $\Psi$, from Assumption H\ref{hypH:coercivity_hyp}. Combining both gives $f(x,\alpha)\leq \liminf_{n} f(x^n,\alpha^n)$.

    This directly implies that $J(t,\lambda;\Pb)=\int f d\Pb$ is lower semicontinuous since the integral of any l.s.c.\ function with at most coercive growth is l.s.c.\ \citep[see, \eg, Lemma 4.3,][]{villani2008optimal}.
\end{proof}
\paragraph{Compactness property.}
To establish the \emph{compactness of the set of $\epsilon$–optimal controls}, two key ingredients are needed: (i) a \emph{metrization of the weak topology in the space of measures}, and (ii) a \emph{representation of the relaxed control problem in terms of SDEs}. These tools allow us to apply compactness criteria for semimartingales, in particular the \emph{Aldous–Rebolledo criterion}, which ensures tightness. Combined with Prokhorov’s theorem, this yields compactness, as shown in the following proposition. The metrization of the weak topology is developed in \Cref{Subsection:weak_topology_measures}, while the SDE representation of relaxed controls is presented in \Cref{Subsection:SDE_representation_martingale_problem}.



\begin{proposition}\label{Prop:Rc_compact}
    Given $(t,\lambda) \in[0,T]\times \Nc(\R^d)$ and $\eps>0$, $\Rc^\eps_{(t,\lambda)}$ is compact in $\Pc^1(\Db^{d}\times \Ac^{\Leb})$.
\end{proposition}
\begin{proof}
    The proof of this lemma breaks into four steps.

    \textit{Step 1.} First, we aim at proving that $\{ \restr{\Pb}{\Db^d}: \Pb \in \CtrlRuleEps_{(t,\lambda)}\}\subseteq\Pc(\Db^b)$ is tight.
    To this end, we verify the Aldous--Rebolledo compactness criterion. For a comprehensive presentation of this criterion, we refer to \citet{rebolledo1980existence}, Theorem 9.4 in \citet{ethier2009markov}, and Theorem 14.11 in \citet{kallenberg1997foundation}. The goal is therefore to show
    \begin{align}
    \label{eq:compact_lemma:limit:aldous_crit}
        \lim_{\delta\downarrow 0} \sup_{\Pb\in \Rc_\lambda} \sup_\tau \E^\Pb\left[\d_{\text{weak*},\R^{d}}(\mu_{(\tau+\delta)\wedge T} ,\mu_{\tau})\right]=0\eqsp,
    \end{align}
    where the innermost supremum is over stopping times $\tau$ valued in $[0,T]$.

    From~\Cref{Prop:representation}, we know there exists an extension $\hat\Omega$ of $\Db^d\times\Ac^{\Leb}$ where $\mu$ can be represented as the solution of~\eqref{SDE:weak}.  This SDE is driven by $\Mc^c$ orthogonal continuous martingale measure on $\hat\Omega\times [0,T]\times\R^d\times A$, with intensity measure $ds\eqsp\mu_s(dx)\eqsp\bar\a_s(x,da)$, and a purely discontinuous martingale measure $\Mc^d$ on $\hat\Omega\times [0,T]\times\R^d\times \R_+ \times A$, with dual predictable projection measure $ds\eqsp\mu_s(dx)\eqsp dz\eqsp\bar\a_s(x,da)$. Applying~\eqref{SDE:weak} to $\varphi_k\in\mathscr{C}_{\R^d}$, we get
    \begin{align*}
        &\langle \varphi_k,\mu_{(s+\delta)\wedge T}\rangle
        \\
        &= \langle \varphi_k,\mu_{s}\rangle +\int_s^{(s+\delta)\wedge T}\int_{\R^d\times A} \big(L \varphi_k (x,\mu,a_r)+
        \gamma(x,\mu,a_r)\left(\partial_s \Phi(1,x,\mu,a_r) - 1\right)\varphi_k(x)\big)\bar\a_r(x,da)\mu_r(dx)dr 
        \\
        &
        ~~+ \int_s^{(s+\delta)\wedge T}\int_{\R^d\times A}
        D\varphi_k(x)\sigma(x,X_r,a) \Mc^{c}(dr,dx,da)
        +  \int_s^{(s+\delta)\wedge T}\int_{\R^d\times \R_+\times A} \sum_{k\geq 0}
        \langle \varphi_k, (k-1)\delta_x\rangle\1_{I_k\left(x ,\mu_{r},a\right)}(z)
        \Mc^d(dr,dx,dz,da)\eqsp.
    \end{align*}
    for $s\in[0,T]$, $k\in\N$. Therefore, to bound the quantity $\E^\Pb\left[|\langle \varphi_k,\mu_{(s+\delta)\wedge T}\rangle - \langle \varphi_k,\mu_{s}\rangle|\right]$, it suffices to bound the last three terms in the r.h.s.\ There is a constant $C>0$ that depends only on $b$, $\sigma$, $\gamma$ and $\Phi$ (which may change from line to line) such that
    \begin{align*}
        &\E^\Pb\Bigg[\Bigg|
        \int_s^{(s+\delta)\wedge T}\int_{\R^d\times A} \big(L \varphi_k (x,\mu,a_r)+ \gamma(x,\mu,a_r)\left(\partial_s \Phi(1,x,\mu,a_r) - 1\right)\varphi_k(x)\big)\bar\a_r(x,da)\mu_r(dx)dr
        \Bigg|\Bigg]
        \\
        & \leq C q_k\E^\Pb\left[ \int_s^{(s+\delta)\wedge T}\left(\langle 1, \mu_u \rangle + \langle |\cdot|, \mu_u \rangle \right)du + \int_s^{(s+\delta)\wedge T}\int_{\R^d\times A}|a| \bar\a_u(x,da)\mu_u(dx)du
        \right]\eqsp.
    \end{align*}
    Applying Burkholder--Davis--Gundy's inequality \citep[see, \eg, Theorem 92,][]{Dellacherie:Meyer:B}, we obtain
    \begin{align*}
        \E^\Pb\Bigg[\Bigg|
        \int_s^{(s+\delta)\wedge T}\int_{\R^d\times A}
                D\varphi_k(x)\sigma(x,X_r,a) \Mc^{c}(dr,dx,da)
        \Bigg|\Bigg]
        \leq
        C q_k\E^\Pb\left[ \int_s^{(s+\delta)\wedge T} \langle 1, \mu_u\rangle du
        \right]
        \eqsp.
    \end{align*}
    Finally, since $\varphi_k\geq 0$, we have
    \begin{align*}
        &\E^\Pb\Bigg[\Bigg|
        \int_s^{(s+\delta)\wedge T}\int_{\R^d\times \R_+\times A} \sum_{k\geq 0}
            \langle \varphi_k, (k-1)\delta_x\rangle\1_{I_k\left(x ,\mu_{r},a\right)}(z)\Mc^d(dr,dx,dz,da)
        \Bigg|\Bigg]\\
        &\leq\E^\Pb\Bigg[\Bigg|
        \int_s^{(s+\delta)\wedge T}\int_{\R^d\times A} \varphi_k(x)\sum_{k\geq 1}
            (k-1)\gamma\left(x ,\mu_{r},a\right)p_k\left(x ,\mu_{r},a\right)\bar\a_r(x,da)\mu_r(dx)dr
        \Bigg|\Bigg]
        \leq C q_k\E^\Pb\left[ \int_s^{(s+\delta)\wedge T} \langle 1, \mu_u\rangle du
        \right]\eqsp.
    \end{align*}
    Combining these inequalities, we get
    \begin{align}
    \label{eq:bound_martingale_terms}
        \begin{split}
            \E^\Pb\left[|\langle \varphi_k,\mu_{(s+\delta)\wedge T}\rangle - \langle \varphi_k,\mu_{s}\rangle|\right]
            &\leq C q_k\E^\Pb\left[ \int_s^{(s+\delta)\wedge T}\left(\langle 1, \mu_u \rangle + \langle |\cdot|, \mu_u \rangle \right)du + \int_s^{(s+\delta)\wedge T}\int_{\R^d\times A}|a| \bar\a_u(x,da)\mu_u(dx)du
            \right]\\
            &\leq
            \delta C q_k\left(\E^\Pb\left[ \sup_{u\in [0,T]}\left(\langle 1, \mu_u \rangle + \langle |\cdot|, \mu_u \rangle \right)\right] + \E^\Pb\left[ \int_s^{(s+\delta)\wedge T}\int_{\R^d\times A}|a| \bar\a_u(x,da)\mu_u(dx)du
            \right]\right)
            \eqsp.
        \end{split}
    \end{align}

    Replacing Itô's formula with \eqref{MartPb:diff-F_f}, and adapting the arguments of Proposition 2.2 and Lemma 2.3 in \citet{ocello2025controlled1}, we obtain that there exists a constant $\bar C>0$, depending only on $T$ and the coefficients $b$, $\sigma$, $\gamma$, and $(p_k)_{k\geq0}$, such that for any $h>0$,
    \begin{align}
        \label{eq:non-explosion-moment1_mass-relaxed}
        \E^{\Pb}\left[\sup_{u\in[t,t+h]}\langle1,\mu_u\rangle\right]\leq& \langle 1,\lambda\rangle e^{C_\gamma C^1_\Phi h}\eqsp,
        \\
        \label{eq:non-explosion-moment2_mass-relaxed}
        \E^{\Pb}\left[\sup_{u\in[t,t+h]}\langle1,
        \mu_u\rangle^2\right]\leq& \langle 1,\lambda\rangle e^{C_\gamma (C^1_\Phi+C^2_\Phi) h}\eqsp,
        \\
        \label{eq:non-explosion-moment2_population-relaxed}
        \E^{\Pb}\left[
            \sup_{u\in[t,t+h]}\langle|\cdot|^2,\mu_u\rangle
        \right]
        \leq& 
        \bar C \Bigg(
            \langle|\cdot|^2,\lambda\rangle
            +\E^{\Pb}\left[ \int_t^{t+h} \langle1,\mu_u\rangle du \right]
            +\E^{\Pb}\left[ \int_t^{t+h}\int_{\R^d\times A} |a|^2 \bar\a_u(x,da)\mu_u(dx)du \right] \bigg)\eqsp,
    \end{align}
    where $|\cdot|^2$ denote the function $x\mapsto |x|^2$. Moreover, adapting Proposition 2.4 in \citet{ocello2025controlled1} to the relaxed setting, we have that
    \begin{align}
        \label{eq:bound_eps_opt_control-relaxed}
        \sup_{\Pb\in \CtrlRuleEps_{(t,\lambda)}} \E^{\Pb}\left[ \int_t^{t+h}\int_{\R^d\times A} |a|^2 \bar\a_u(x,da)\mu_u(dx)du \right]<\infty\eqsp.
    \end{align}
    Therefore, from the previous bounds and \eqref{eq:bound_martingale_terms}, we obtain $\E^\Pb \big[ |\langle \varphi_k,\mu_{(s+\delta)\wedge T}\rangle - \langle \varphi_k,\mu_{s}\rangle|\big]\leq C q_k\delta$.
    Multiplying for $\frac{1}{2^k q_k}$, summing over $k\in\N$ and applying the monotone convergence theorem, we get $\E^\Pb\left[\d_{\R^d}(\mu_{(s+\delta)\wedge T} ,\mu_{s})\right]\leq\delta C$, which gives us~\eqref{eq:compact_lemma:limit:aldous_crit}.

    \textit{Step 2.} Secondly, we prove that $\left\{ \restr{\Pb}{\Db^d}: \Pb \in \CtrlRuleEps_{(t,\lambda)}\right\}\subseteq\Pc^1(\Db^b)$ is relatively compact. Combining the bound~\eqref{eq:non-explosion-moment1_mass-relaxed} with \eqref{eq:non-explosion-moment2_population-relaxed} and \eqref{eq:bound_eps_opt_control-relaxed}, we get
    \begin{align*}
        \sup_{\Pb\in \Rc^{\eps}_{(t,\lambda)}} \E^{\Pb}\left[ \sup_{u\in[t,T]}\int_{\R^d} |x|^2 \mu_u(dx) \right]<\infty\eqsp.
    \end{align*}
    This bound, together with~\eqref{eq:non-explosion-moment2_mass-relaxed} and~\eqref{eq:bound_d_p_E}, gives that
    \begin{align}
    \label{eq:bound_eps_opt_traj-relaxed}
        \sup_{\Pb\in \Rc^{\eps}_{(t,\lambda)}} \E^{\Pb}\left[ \sup_{u\in[t,T]} \d^2_{2,\R^d}(\mu_u,\delta_{0})\right]<\infty\eqsp.
    \end{align}
    Putting together \textit{Step 1} and this bound, we have from Corollary B.2 of \citet{Lacker:MFgames:2015} that $\left\{ \restr{\Pb}{\Db^d}: \Pb \in \CtrlRuleEps_{(t,\lambda)}\right\}\subseteq\Pc^1(\Db^b)$ is relatively compact.

    \textit{Step 3.} From the first step, we have that $\left\{ \Pb\circ\mu^{-1}: \Pb \in \CtrlRuleEps_{(t,\lambda)}\right\}$ is tight in $\Pc^1(\Db^b)$. Adding this to~\eqref{eq:bound_eps_opt_control-relaxed} and~\eqref{eq:bound_eps_opt_traj-relaxed}, we have that $\left\{\restr{\Pb}{\Ac^{\Leb}}:\Pb\in\CtrlRuleEps_{(t,\lambda)}\right\}$ is compact in $\Pc^1(\Ac^{\Leb})$. This entails that $\CtrlRuleEps_{(t,\lambda)}$ is relatively compact in $\Pc^1(\Db^{d}\times \Ac^{\Leb})$ since $\Big\{ \restr{\Pb}{\Db^d}: \Pb \in \CtrlRuleEps_{(t,\lambda)}\Big\}$ and $\Big\{\restr{\Pb}{\Ac^{\Leb}}:\Pb\in\CtrlRuleEps_{(t,\lambda)}\Big\}$ are relatively compact in $\Pc^1(\Db^d)$ and $\Pc^1(\Ac^{\Leb})$ respectively.

    \textit{Step 4.} Finally, we prove $\CtrlRuleEps_{(t,\lambda)}$ is closed. To do that, we show that $\Pb^\infty$ belongs to $\CtrlRuleEps_{(t,\lambda)}$, for $\Pb^n \to \Pb^\infty$ in $\Pc^1(\Db^{d} \times \Ac^{\Leb})$, with $\Pb^n\in \CtrlRuleEps_{(t,\lambda)}$. Since $\mu_t=\lambda$ under $\Pb^n$, the same is true under $\Pb^\infty$. Analogously, since $\Pb^n(\alpha\in \Ac^{\Leb, \mu})=1$, the same is true under $\Pb^\infty$. For any $F\in C^2_b(\R)$ and $\varphi\in C^2_b(\R^d)$ and $\Pb\in \Pc^1(\Db^{d} \times \Ac^{\Leb})$, define
    $M^{\Pb,F_\varphi}_s : \Db^{d} \times \Ac^{\Leb} \to \R$ by
    \begin{align*}
        M^{\Pb,F_\varphi}_s(\x,\alpha) =& F_\varphi(\x_s) - \int_t^s \int_{\R^d\times A} \Lc F_\varphi(y,\y_u, a)\bar\alpha_u(y,da)\y_u(dy)\delta_{\y_u=\x_u}du\eqsp.
    \end{align*}
    Recalling the definition of $\Lc$, we see that there exists a constant $C>0$ depending only on the bounds of $F$, $\varphi$ and the constants $C_b$, $C_\sigma$, $C_\gamma$ such that 
    \begin{align*}
        \left|\Lc F_\varphi(y,\lambda, a)\right|\leq C(1 + |x|+|a|)\eqsp.
    \end{align*}
    This implies
    \begin{align*}
        \left|M^{\Pb,F_\varphi}_s(\x,\alpha)\right|\leq C\left(
            1 + \sup_{u\in[t,T]}
            \d_{1,\R^d}\left(\x_u, \delta_0\right) +\int_t^T \int_{\R^d\times A}|a|\bar\alpha_u(x,da)\x_u(dx)du
        \right)\eqsp,
    \end{align*}
    which is uniformly bounded for $\Pb\in \CtrlRuleEps_{(t,\lambda)}$ from~\eqref{eq:bound_eps_opt_control-relaxed} and~\eqref{eq:bound_eps_opt_traj-relaxed}.
    Combining this with the continuity of $b$, $\sigma$, $\gamma$ and $p_k$, for $k\in\N$,
    since $\Pb^n\to \Pb^\infty$ in $\Pc^1(\Db^{d} \times \Ac^{\Leb})$, it follows that
    \begin{align*}
        \E^{\Pb^\infty}\left[\left(M^{\Pb^\infty,F_\varphi}_{s+u} -M^{\Pb^\infty,F_\varphi}_s\right)\mathfrak{f} \right]=\lim_{n\to \infty}\E^{\Pb^n}\left[\left(M^{\Pb^n,F_\varphi}_{s+u} -M^{\Pb^n,F_\varphi}_s\right)\mathfrak{f} \right]\eqsp,
    \end{align*}
    for every $s\in[t,T]$, $u\geq0$ such that $s+u\leq T$, any $F\in C^2_b(\R)$ and $\varphi\in C^2_b(\R^d)$, and any bounded continuous function $\mathfrak{f}$ on $\Db^{d} \times \Ac^{\Leb}$, measurable with respect to $\sigma\left(\mu_u,\bar\a_u : u \in [t,s]\right)$. Since $\Pb^n\in \CtrlRuleEps_{(t,\lambda)}$, the process
    $\left(M^{\Pb^n,F_\varphi}_s(\mu,\a)\right)_{s\in [0,T]}$ is a martingale under $\Pb^n$, and the above quantity is zero. This shows that
    $\left(M^{\Pb^\infty, \varphi}_s(\mu,\a)\right)_{s\in [0,T]}$ is a martingale under $\Pb^\infty$, and so $\Pb^\infty\in \Rc_{(t,\lambda)}$.

    Moreover, by~\Cref{Lemma:J_continuous} we get since $J$ is lower semicontinuous. Therefore,
    \begin{align*}
        J\left(t,\lambda;\Pb^\infty\right)\leq \liminf_{n\to\infty}J(t,\lambda;\Pb^n)\leq \valueRelaxed(t,\lambda) +\eps\eqsp,
    \end{align*}
    which means that $\Pb^\infty\in \CtrlRuleEps_{(t,\lambda)}$.
\end{proof}

\paragraph{Existence of optimal controls.}
We can now establish the existence of an optimal control, which is the main result of this section
\begin{theorem}\label{thm:existence-opt-controls}
    For $(t,\lambda) \in[0,T]\times \Nc(\R^d)$, there exists an optimal control $\beta^*\in \Rc^\mathfrak{s}_{(t,\lambda)}$ such that
    \begin{align}\label{eq:min_atteint}
        \valueStrong(t,\lambda)=J(t,\lambda;\beta^*)\eqsp.
    \end{align}
\end{theorem}
\begin{proof}
    Fix $\eps>0$. We have that $\inf_{\Pb\in \Rc_{(t,\lambda)}}J(t,\lambda;\Pb)=\inf_{\Pb\in \Rc^\eps_{(t,\lambda)}}J(t,\lambda;\Pb)$. By~\Cref{Prop:Rc_compact}, $\Rc^\eps_{(t,\lambda)}$ is compact and, by~\Cref{Lemma:J_continuous}, $J$ is lower-semicontinuous. Therefore, since $\valueStrong(t,\lambda)=\valueRelaxed(t,\lambda)$ is the supremum of a continuous function over a nonempty compact set, it exists $\Pb^*\in\Rc_{(t,\lambda)}$ such that $\valueStrong(t,\lambda)=J(t,\lambda;\Pb^*)$. From~\Cref{Lemma:F_X-measurability} and~\Cref{Prop:atom_ctrlVSnat_ctrl}, under Assumption H\ref{Assumption:Convex}, we have the existence of optimal weak control $\mathfrak{a}^*$ such that $J(t,\lambda;\mathfrak{a}^*)\leq J(t,\lambda;\Pb^*)$. Immerging this weak control in the class of strong controls, as done in the proof of \Cref{thm:equivalence_forms}, we find $\beta^*$ that satisfies~\eqref{eq:min_atteint}.
\end{proof}

\section{Conclusion}\label{Section:Conclusion}

This article is the companion work to \citet{ocello2025controlled1}, where we began the analysis of controlled branching diffusions in the mean-field regime. Here, we complete that study by addressing the existence of optimal controls, under a relaxed formulation point of view.

Under the mean-field interaction assumption, we formulated the control problem through its associated martingale problem. By developing the notions of natural and weak controls, we reduced the scope of the problem while preserving generality. Through a Filippov-type convexity condition, we proved the equivalence between relaxed, weak, and strong controls. Shifting to the framework of control rules, we showed that the optimization problem can be confined to a compact set and that the cost functional is lower semicontinuous. This guarantees the existence of an optimal solution for both the relaxed and strong formulations.

Beyond existence, our analysis sets the stage for further developments. The relaxed formulation provides a rigorous foundation for numerical methods, such as Markov chain approximations, finite-difference schemes for HJB equations, and reinforcement learning algorithms. Moreover, it opens the way to scaling limits of controlled populations: in particular, \citet{ocello2025controlled} studies the emergence of controlled superprocesses as natural limits. Combining this approach with recent advances on heterogeneous systems \citep{lacker2023label,de2024mean,coppini2025nonlinear,de2025linear} would allow the study of large-scale interacting populations with non-uniform structures, offering a rich and promising direction for future research.

\textbf{Acknowledgements.} This work is supported by Hi! PARIS and ANR/France 2030 program (ANR-23-IACL-0005). I gratefully acknowledge my PhD supervisor Idris Kharroubi for supervising this work.

\appendix

\section{Appendix}

\subsection{Weak topology and metrization on spaces of measures}
\label{Subsection:weak_topology_measures}

Since $\d_{1,\Ec}$ is a Wasserstein type distance and we have the bound~\eqref{eq:bound_d_p_E}, the results from Appendix B of \citet{Lacker:MFgames:2015} can be naturally extended to this setting. As the primary focus is on convergence in weak* topology in the first part, we will examine an alternative metrization that is simpler than $\d_{1,\Ec}$.

A family $\mathscr{C}\subseteq C_b(\Ec)$ is said to be \textit{separating} if, whenever $\langle \varphi,\lambda\rangle=\langle \varphi,\lambda'\rangle$, for all $\varphi\in \mathscr{C}$, and some $\lambda, \lambda'\in M(\Ec)$, we necessarily have $\lambda=\lambda'$. Since $\Ec$ is Polish, from the Portmanteau theorem \citep[see, \eg, Theorem 1.1.1,][]{SV97}, the set of uniformly continuous functions, for metric equivalent to $d$, is separating. Using Tychonoff's embedding theorem \citep[see, \eg, Theorem 17.8,][]{General:Topology}, $C_b(\Ec)$ is also separable.
Therefore, there exists a countable and separating family $\mathscr{C}_\Ec = \left\{\varphi_k, k \in\N\right\}$ subset of $C_b(\Ec)$ such that the function $\Ec\ni x\mapsto 1$ belongs to $\mathscr{C}_\Ec$ and $\|\varphi_k\|_\infty := \sup_\Ec|\varphi_k| \leq 1$, for all $k \in\N$ since multiplying by a positive constant do not impact the property of being separating. With the use of this family, 
\begin{align*}
    \d_{\text{weak*},\Ec}(\lambda,\lambda') = \sum_{\varphi_k\in\mathscr{C}_\Ec}\frac{1}{2^k}\left|\langle \varphi_k,\lambda\rangle-\langle \varphi_k,\lambda'\rangle\right|\eqsp,
\end{align*}
for $\lambda,\lambda'\in M(\Ec)$. As in Theorem 1.1.2 of \citet{SV97}, this distance $\d_{\text{weak},\Ec}$ induces on $M(\Ec)$ the weak* topology. Whenever $\Ec=\R^d$, we adjust this metric to take into account useful differential properties. Let $\mathscr{C}_{\R^d}$ be taken as a subset of $C^2_b(\R^d)$, the set of real functions with bounded, continuous derivatives over $\R^d$ up to order two. Without loss of generality, since $C^2$ is dense in $C^0$, we suppose this set to be separating under local uniform convergence \citep[application of Theorem 8.14 in ][]{folland1999real}. Moreover, since $\x\mapsto1$ belongs to $\mathscr{C}_{\R^d}$, adding a constant or multiplying by a non-negative constant to each function does not change the property of being a separating set, we assume $\varphi_k\geq 0$. We define the distance
\begin{align*}
    \d_{\text{weak*},\R^d}(\lambda,\lambda') = \sum_{\varphi_k\in\mathscr{C}_{\R^d}}\frac{1}{2^k q_k}\left|\langle \varphi_k,\lambda\rangle-\langle \varphi_k,\lambda'\rangle\right|\eqsp,
\end{align*}
with $q_k=\max\{1,||D \varphi_k||_\infty, ||D^2\varphi_k||_\infty\}$.

\subsection{SDE representation of the martingale problem}
\label{Subsection:SDE_representation_martingale_problem}

We now represent relaxed controls through the SDE point of view. This description is particularly useful for establishing the existence of optimal controls. The construction relies on \emph{martingale measures}, defined on an extended probability space; we briefly recall their definition here, referring to \citet{ElK_Meleard:mart_measure} and the monograph of \citet{Walsh} for further details.

\begin{definition}
Let $(G,\Gc)$ be a Lusin space with its $\sigma$-algebra, and $(
\Omega, \mathscr{F}, \Pb,\mathfrak{F} =\{\mathscr{F}_s\}_{s\geq0})$ a filtered space satisfying the usual condition, where we define $\Pc$ the predictable $\sigma$-field. A process $\Mc$ on $\Omega\times[0,T]\times \Gc$ is called  \emph{martingale measure on $G$} if
\begin{enumerate}[(i)]
\item $\Mc_0(E)=0$ a.s. for $E\in\Gc$;
\item $\Mc_t$ is a $\sigma$-finite, $L^2(\Omega)$-valued measure for all $t\in [0,T]$;
\item $\left(\Mc_t(E)\right)_{t\in[0,T]}$ is an $\mathfrak{F}$-martingale for $E\in\Gc$.
\end{enumerate}
We say that $\Mc$ is \emph{orthogonal} if the product $\Mc_t(E)\Mc_t(E')$ is a martingale for two disjoint sets $E,E'\in \Gc$. We also say, on one hand, that is \emph{continuous} if $\left(\Mc_t(E)\right)_{t\geq0}$ is continuous, \emph{purely discontinuous}, on the other hand, if $\left(\Mc_t(E)\right)_{t\geq0}$ is a purely discontinuous martingale for $E\in \Gc$.
\end{definition}

For a strong representation of relaxed controls, we rely on the notion of \emph{predictable projection} and \emph{intensity} that we briefly recall.
For an $\R$-valued $\mathfrak{F}$-adapted process $Y$, there exists \citep[see, \eg, Theorem 2.28, Chapter I,][]{jacod2013limit} a $(-\infty,\infty]$-valued process, called the \emph{predictable projection} of $Y$ and denoted by ${}^P Y$. It is determined uniquely up to a negligible set by the following two conditions:
\begin{enumerate}[(i)]
\item it is predictable;
\item ${}^P Y_T = \E^{\Pb}\left[Y_T|\mathscr{F}_{T-}\right]$ on $\{T < \infty\}$, for all predictable stopping times $T$.
\end{enumerate}

For a continuous orthogonal martingale measure $\Mc$ on $G$, there exists a random, predictable real-valued measure $I$ on $\Bc([0,T])\otimes \Gc$, called \emph{intensity} of $\Mc$, defined by: $\left[\Mc(E)\right]_s= \int_0^t\int_E I(dx, ds)$ $\Pb$--a.s., for all $t > 0$. We can construct a stochastic integral with respect to $\Mc$, for all functions $\varphi$ defined on $\Omega\times[0,T]\times G,$ $\Pc\otimes \Gc$ measurable, such that
\begin{align*}
    \E^{\Pb}\left[\int_0^t\int_E \varphi^2(\omega,s,x) I(\omega,dx, ds)\right]<\infty\eqsp,
\end{align*}
denoted by $\int_0^t\int_E \varphi(s,x) \Mc(dx, ds)$. We refer to Chapter 2 of \citet{Walsh} for the proofs.

The representation of these processes is grounded in the representation theorems for continuous and purely discontinuous martingale measures, as done in \citet{Roelly-Meleard:discont_measure_val_branching}. We apply her construction in our context and get the following proposition.

\begin{proposition}
\label{Prop:representation}
    Let $\Cc = ( \Omega, \mathscr{F}, \Pb,\mathfrak{F} =\{\mathscr{F}_s\}_{s\geq0}, X, \alpha ) \in \RelCtrl_{(t,\lambda)}$. There exists an extension \begin{align*}
        \left(
    \hat\Omega=\Omega\times\tilde\Omega, \hat{\mathscr{F}}=\mathscr{F}\otimes\tilde{\mathscr{F}}, \hat\Pb=\Pb\otimes\tilde\Pb,\left\{\hat{\mathscr{F}}_s=\mathscr{F}_s\otimes\tilde{\mathscr{F}}_s\right\}_s\right)
    \end{align*}
    of $(\Omega, \mathscr{F}, \Pb,\mathfrak{F})$, where we naturally extend $X$ and $\alpha$, that satisfies the following properties.
    \begin{enumerate}
        \item $(\hat\Omega,\hat{\mathscr{F}},\hat{\mathfrak{F}},\hat\Pb)$ is a filtered probability space supporting two orthogonal $\hat{\mathfrak{F}}$-martingale measures: a continuous one, $\Mc^c$, on $\hat\Omega\times [0,T]\times\R^d\times A$, with intensity measure $dsX_s(dx)\bar\alpha_s(x,da)$, and a purely discontinuous one, $\Mc^d$, on $\hat\Omega\times [0,T]\times\R^d\times \R_+ \times A$, with dual predictable projection measure $dsX_s(dx)dz\bar\alpha_s(x,da)$.
        \item $\hat\Pb\circ X_t^{-1}=\lambda$.
        \item $\hat\Pb(\alpha\in \Ac^{\Leb, X})=1$.
        \item $X$ satisfies the following dynamics
        \begin{align}
        \label{SDE:weak}
            \begin{split}
                \langle f, X_s\rangle = \langle f, \lambda\rangle &+
                \int_t^s\int_{\R^d\times A}   \big(L f (x,X_r,a)+\gamma(x,X_r,a)\left(\partial_s \Phi(1,x,X_r,a) - 1\right)f(x)\big)\bar\alpha_r(x,da)X_r(dx)dr
                \\
                &+ \int_t^s\int_{\R^d\times A}
                Df(x)\sigma(x,X_s,a) \Mc^c(dr,dx,da)
                \\
                &+ \int_t^s\int_{\R^d\times \R_+\times A} \sum_{k\geq 0}
                \langle f, (k-1)\delta_x\rangle\1_{I_k\left(x ,X_{r},a\right)}(z)\Mc^d(dr,dx,dz,da)
                \eqsp,
            \end{split}
        \end{align}
        for all $f\in C^\infty_b(\R^d)$ and all $[t,s]\subseteq[0,T]$.
    \end{enumerate}
\end{proposition}

\begin{proof}
    We follow the ideas in Theorem 2.7 of \citet{Roelly-Meleard:discont_measure_val_branching} and Theorem 2.9 of \citet{Roelly-Meleard:discont_measure_val_branching} to characterize the martingale $\bar M^{f}_s$ in~\eqref{MartPb:diff-var_finie}. As proven in Theorem 4.18 of \citet{jacod2013limit}, every square integrable martingale starting at $0$ can be uniquely decomposed in the sum of a continuous martingale $\bar M^{f,c}$ and a purely discontinuous martingale $\bar M^{f,d}$, which is the compensated sum of its jumps. We show the connection of these two processes with $X$ and $\alpha$.

    First, we focus on $\bar M^{f,d}$. Since a purely discontinuous martingale $\bar M^{f,d}$ is the compensated sum of its jumps, we look at $\Delta X_s = X_s-X_{s-}$. Let $\tilde N$ be the Lévy system of $X$, \ie, a measure on $M^1(\R^d)\times \R_+$ given by $N_s(X_s,dv)ds$ where $N_s(\bar X,dv)$ is the image measure of the measure $\nu_s(x,\bar X,du)\bar X(dx)$ by the mapping $(u,x)\mapsto u \delta_x$ from $\R_+\times  \R^d$ to $M^1(\R^d)$, and a certain kernel $\nu$. Comparing the last term in expressions~\eqref{MartPb:diff-F_f} and Théorème 7 (4) of \citet{ElK:Roelly:Levy-Khintchine_repr}, we identify $\nu$ as
    \begin{align*}
        \nu_s(x,\lambda,dz) = \int_A \sum_{k\geq 0}\left(k-1\right) \1_{I_k(x,\lambda,a)}(z) \bar \alpha_s(x,da)dz\eqsp.
    \end{align*}
    This means that, for $F$ bounded positive measurable function on $\R_+\times M^1(\R^d)$, we have that
    \begin{align}
        &\sum_{t< r\leq s}F( r,\Delta X_{ r}) \1_{\{\Delta X_{r} \neq 0\}} 
        - \int_t^s \int_{\R^d}   \int_{(0,\infty)} \int_A \sum_{k\geq 0} F( r, (k-1) \delta_x)  \1_{I_k(x,X_{r},a)}(z) \bar \alpha_r(x,da) dz X_r(dx) dr
        \\
        &
        =\sum_{t< r\leq s}F( r,\Delta X_{r}) \1_{\{\Delta X_{r} \neq 0\}} \label{SDE:repr_Levy_system}
        -\int_t^s \int_{\R^d\times A}   \sum_{k \geq 0} F( r, (k-1) \delta_x) \gamma(x,X_r,a)p_k(x,X_r,a)  \bar\alpha_r(x,da)X_r(dx) dr 
    \end{align}
    is a $\F$-martingale. With this description of $\nu$ and $N_s(X_s,dv)ds$, we use Proposition 2.8 of \citet{Roelly-Meleard:discont_measure_val_branching} to prove that we satisfy the hypothesis of Theorem 12 of \citet{ElK:Lepeltier:proc_ponct}. Therefore, there exists an extension $\big(
    \bar\Omega^1=\Omega\times \Omega^1, \bar\Fc^1=\Fc\otimes \Fc^1, \bar\Pb^1=\Pb\otimes \Pb^1, \left\{\bar\Fc_s^1=\Fc_s\otimes \Fc_s^1\right\}_{s\geq0}\big)$, and martingale measures $\Mc^d$ on $[0,T]\times \R^d\times\R_+\times A$ in it, such that its dual predictable projection measure is $dr\eqsp X_r(dx)dz\eqsp\bar\alpha_r(x,da)$, and 
    \begin{align*}
        \bar M^{f,d}_s = \int_t^s\int_{\R^d\times \R_+\times A} \sum_{k\geq 0}
        \langle f, (k-1)\delta_x\rangle\1_{I_k\left(x ,X_{r},a\right)}(z)\Mc^d(dr,dx,dz,da)\eqsp.
    \end{align*}

    We focus now on $\bar M^{f,c}$. The first term in~\eqref{MartPb:diff-quadratic_var} comes from the continuous martingale, \ie,
    \begin{align*}
        \left[\bar M^{f,c}\right]_s = \int_t^s \int_{\R^d\times A}    \text{Tr}\left(\sigma\sigma^\top(x,X_r,a)D\varphi D\varphi^\top(x)\right) \bar\alpha_r(x,da)X_r(dx) dr
        \eqsp.
    \end{align*}
    Since $\sigma\in L^2( X_s(dx) \alpha_s(da)ds)$, from Theorem III-7 of \citet{ElK_Meleard:mart_measure}, there exist an extension $\big(
    \bar\Omega^2=\bar\Omega^1\times \Omega^2, \bar\Fc^2=\bar\Fc^1\otimes \Fc^2, \bar\Pb^2=\bar\Pb^1\otimes \Pb^2, \left\{\bar\Fc_s^2=\bar\Fc_s^1\otimes \Fc_s^2\right\}_{s\geq0}\big)$, and a continuous martingale measure $\Mc^c$ on $[0,T]\times\R^d\times A$ on this space, such that its intensity is $ds X_s(dx)\bar\alpha_s(x,da)$, and we have
    \begin{align*}
        \bar M^{f,c}_s = \int_t^s\int_{\R^d\times A}  Df(x)\sigma(x,X_r,a)\Mc^c(dr,dx,da)\eqsp.
    \end{align*}
    The imposed dependence on $X$ and $\alpha$ over $\Mc^d$ and $\Mc^c$ implies that~\eqref{SDE:weak} is satisfied.

    Conversely, if a $M^1(\R^d)$-valued process satisfies~\eqref{SDE:weak}, applying Itô's formula, we have~\eqref{MartPb:diff-exp}.
\end{proof}

\bibliographystyle{plainnat}
\bibliography{main}

\end{document}